\title{ ~~\\ Sequences of enumerative geometry: congruences
  and asymptotics}
\author{Daniel B. Gr\"unberg, Pieter Moree\\ \small{Appendix by Don Zagier}}
\documentclass[12pt]{article}
\usepackage{amssymb, latexsym, amsfonts}
\def\@ptsize{2}
\newtheorem{Thm}{Theorem}

\newcommand{\half}{{1\over 2}}
\newtheorem{lem}{Lemma}
\newtheorem{cor}{Corollary}

\newcommand{\p}{\partial}

\newtheorem{prop}{Proposition}
\newcommand{\qed}{\hfill $\Box$}
\newcommand{\C}{\mathbb{C}}
\newcommand{\N}{\mathbb{N}}
\renewcommand{\P}{\mathbb{P}}
\newcommand{\Z}{\mathbb{Z}}

\begin{document}
\date{}
\maketitle
{\def\thefootnote{}
\footnote{{\it Mathematics Subject Classification (2000)}. 14N10,
11A07, 41A60}}
\begin{abstract}
\noindent 
We study the integer sequence $v_n$ of numbers of lines in hypersurfaces
of degree $2n-3$ of $\P^n$, $n>1$.  
We prove a number of congruence properties of these numbers of several different types. 
Furthermore, the asymptotics of the $v_n$ are described (in an appendix by Don Zagier). 
An attempt is made at a similar analysis of two other enumerative sequences: the numbers of
rational plane curves and the numbers of instantons in the quintic
threefold.
\end{abstract}
\noindent We study the sequence of numbers of lines in a hypersurface of degree
$D=2n-3$ of $\P^n$, $n>1$.  The sequence is defined by (see e.g.
\cite{F-84})
\begin{equation}
\label{een}
v_n:=\int_{G(2,n+1)} c_{2n-2}({\rm Sym}^D Q),
\end{equation}
where $G(2,n+1)$ is the Grassmannian of $\C^2$
subspaces of $\C^{n+1}$ (i.e. projective lines in $\P^n$) of
dimension $2(n+1-2)=2n-2$, $Q$ is the bundle of linear
forms on the line (of rank $r=2$, corresponding to a particular point of the
Grassmannian), and Sym$^D$ is its $D$th symmetric product -- of rank
${D+r-1 \choose r-1} = D-1 =2n-2$. The top Chern class (Euler class)
$c_{2n-2}$ is the class dual to the 0-chain (i.e. points) corresponding
to the zeros of the bundle Sym$^D(Q)$, i.e. to the vanishing of a degree
$D$ equation in $\P^n$; this is the geometric requirement that the
lines lie in a hypersurface.\\
\indent The integral (\ref{een}) can actually be written as a sum:
\begin{equation}
\label{sommie}
v_n=\sum_{0\leq i<j
  \leq n} {\prod_{a=0}^D (a w_i + (D-a) w_j) \over
  \prod_{0\leq k\leq n,~k\neq i,j} (w_i -w_k)(w_j -w_k) },
\end{equation}
where $w_0,\ldots,w_n$ are arbitrary complex variables.
This is a consequence of a localisation formula due to
Atiyah and Bott from equivariant cohomology, which says that only the
(isolated) fixed points of the $(\C^*)^{n+1}$ action contribute to the defining
integral of $v_n$. Hence the sum. For the first few values of $n$ computation
yields $v_2=1$, $v_3=27$, $v_4=2875$, $v_5=698005$, $v_6=305093061$.\\
\indent D. Zagier gave a simple proof that the right hand side of (\ref{sommie}) 
is independent of $w_0,\ldots,w_n$ (as it must be for (\ref{sommie}) to hold), and
that in fact it can be replaced by the much simpler formula
\begin{equation} \label{eq:defn}
v_n = \Big[ (1-x) \prod_{j=0}^{2n-3} (2n-3-j+jx) \Big]_{x^{n-1}}
\end{equation}
where the notation $[\dots]_{x^n}$ means the coefficient of $x^n$. 
In fact, formula (\ref{sommie}) was proved in a very different way using
methods from Schubert calculus by B.L. van der Waerden, who established it in part 2 of his
celebrated 20 part `Zur algebraischen Geometrie' series of papers \cite{vdW0, vdWS}. The number of linear subspaces of
dimension $k$ contained in a generic hypersurface of degree $d$ in $\P^n$, when it is finite, 
can be likewise expressed as the coefficient of a monomial in a certain polynomial
in several variables, see e.g. \cite[Theorem 3.5.18]{Mani}.\\
\indent Zagier also gave the formula
\begin{equation}
\label{vierie}
v_n\sim \sqrt{\frac{27}{\pi}} (2n-3)^{2n-7/2} \Big( 1-\frac{9}{8n}
-\frac{111}{640 n^2} -\frac{9999}{25600 n^3}+\cdots \Big),
\end{equation}
where the right hand side is an asymptotic expansion in powers of $n^{-1}$ with rational
coefficients that can be explicitly computed. The proof of this formula, as well as the derivation
of (\ref{eq:defn}) from (\ref{sommie}), can be found in the appendix.\\
\indent The remaining results, summarized in 
Theorem \ref{congthm} and Theorem \ref{extra1},  are concerned with congruence properties of the
numbers $v_n$. In this context it turns out to be convenient to
define $v_1=1$ (even though there 
is no such thing as a 
hypersurface in $\P^1$ of degree $-1$) and even more remarkably $v_0=-1$. We do not doubt that the 
congruence results presented here form
only the tip of an iceberg.\\
\indent A first version of this paper was single authored by the first author. Indeed,
the present version of this paper is similar to the first one, except
for sections 2 and 3 which have been greatly revised and expanded by the
second author. The conjectures outside these two sections are due to the
first author alone. Sections 4 and 5 were revised by both the second
author and Don Zagier.

\section{Introduction}

The motivating idea behind this paper is the expectation 
that certain problems in enumerative geometry are
coupled to modularity.  This is a recurrent theme in string theory,
where partition functions have often an enumerative interpretation as
counting objects (instantons, etc) and must satisfy the condition of
modularity covariance in order to obtain the same amplitude when two
worldsheets have the same intrinsic geometry.

Modular forms, as is well known, have Fourier coefficients satisfying
many interesting congruences (think of Ramanujan's congruences for
partitions, or for his function $\tau(n)$).  The same can happen for the 
coefficients of expansions related to modular forms, e.g.~the expansions
$y=\sum A_nx^n$ obtained by writing a modular form $y$ (locally) as a
power series in a modular function $x$.  For instance, the famous Ap\'ery
numbers related to Ap\'ery's proof of the irrationality of $\zeta(3)$ are
obtainable in this way \cite{Beu} and satisfy many interesting congruences
\cite{SB}.  The numbers appearing in the context of mirror symmetry,
Picard-Fuchs equations for Calabi-Yau manifolds, Gromov-Witten invariants
and similar problems of enumerative geometry are sometimes related to
modular forms and sometimes not, so we can reasonably hope for interesting
congruence properties in these contexts also.
In Section 2 we shall find astonishingly many congruences for our sequence $v_n$. 
We shall first draw a few
tables for congruences mod 2,3,4,5 or 11, and then summarize the
observed congruences.  In Section 3 we prove those
congruences by elementary means starting from (\ref{eq:defn}), and a few conjectures will 
be formulated. 
Sequences of numbers coming from modular forms also often have interesting
asymptotic properties and we therefore wish to study this, too.  In
Section 4 we find the asymptotic properties of the $v_n$ numerically by
using a clever empirical trick shown to us by Don Zagier which we call
the asymp$_k$ trick.  (A rigorous proof of these asymptotics, as already
mentioned above, was also provided by him and is reproduced in the
appendix.) Section 5 presents congruences and asymptotics for two
further examples of enumerative sequences, without proofs. The sequence of rational curves on the plane and the
sequence of instantons on the quintic threefolds partly partly mimic
the behaviour of the original sequence of lines in hypersurfaces.  

\section{Congruences}

We will consider the sequences $\{v_n\; ({\rm mod~}k)\}_{k=1}^{\infty}$ for
some small values of $k$; that is, we study the reduction of the integers
$v_n$ modulo $k$. It turns out to be instructive to order the $v_n$ modulo $k$ 
in a table. Each table has $k$ columns. The $i$th column $(i=1,\ldots,k)$ gives
the values of $v_{lk+i}\; ({\rm mod~}k)$ $(l\ge 0)$.

For instance, the first few tables at $k=2,3,4,\dots$ look like \\
\underline{$k=2$}:
\begin{tabular}{cccc cccc}
  1 & 1 & 1 & 1 & 1 & 1 & 1 & \dots\\
  1 & 1 & 1 & 1 & 1 & 1 & 1 & \dots
\end{tabular} \hspace{0.5cm}
\underline{$k=3$}:
\begin{tabular}{cccc cccc cc}
  1 & 1 & 2 & 2 & 2 & 0 & 0 & 0 & 2 & \dots\\
  1 & 1 & 2 & 2 & 2 & 0 & 0 & 0 & 2 & \dots\\
  0 & 0 & 0 & 0 & 0 & 0 & 0 & 0 & 0 & \dots
\end{tabular}\\ \\ \\
\underline{$k=4$}:
\begin{tabular}{cccc cccc}
  1 & 1 & 1 & 1 & 1 & 1 & 1 & \dots\\
  1 & 1 & 1 & 1 & 1 & 1 & 1 & \dots\\
  3 & 3 & 3 & 3 & 3 & 3 & 3 & \dots\\
  3 & 3 & 3 & 3 & 3 & 3 & 3 & \dots
\end{tabular} \hspace{0.5cm}
\underline{$k=5$}:
\begin{tabular}{cccc cccc cc}
  1 & 1 & 4 & 4 & 4 & 4 & 4 & 2 & 3 & \dots\\
  1 & 1 & 4 & 4 & 4 & 4 & 4 & 2 & 3 & \dots\\
  2 & 0 & 0 & 2 & 3 & 2 & 1 & 4 & 1 & \dots\\
  0 & 0 & 0 & 0 & 0 & 0 & 0 & 0 & 0 & \dots\\
  0 & 1 & 0 & 0 & 0 & 0 & 4 & 0 & 0 & \dots
\end{tabular}\\
The first table ($k=2$) says that all the $v_n$ are odd integers.
We shall mostly be interested in the tables for prime $k$. Here is a
typical prime table:\\ \\
\underline{$k=11$}: {\small
\begin{tabular}{cccc cccc cccc ccc}
  1 & 1 & 10 & 10 & 10 & 10 & 10 & 10 & 10 & 10 & 10 & 10 & 10 & 8 & \dots\\
  1 & 1 & 10 & 10 & 10 & 10 & 10 & 10 & 10 & 10 & 10 & 10 & 10 & 8 & \dots\\
  5 & 9 & 10 & 7 & 8 & 6 & 10 & 2 & 7 & 8 & 6 & 10 & 8 & 5 & \dots\\
  4 & 1 & 5 & 8 & 6 & 7 & 10 & 8 & 2 & 6 & 7 & 10 & 8 & 8 & \dots\\
  0 & 9 & 3 & 2 & 2 & 0 & 0 & 0 & 0 & 0 & 0 & 0 & 1 & 5 & \dots\\
  9 & 3 & 10 & 0 & 1 & 4 & 8 & 10 & 7 & 6 & 2 & 8 & 10 & 7 & \dots\\
  0 & 0 & 0 & 0 & 0 & 0 & 0 & 0 & 0 & 0 & 0 & 0 & 0 & 0 & \dots\\
  0 & 2 & 5 & 1 & 1 & 0 & 0 & 0 & 0 & 0 & 0 & 0 & 8 & 10 & \dots\\
  0 & 2 & 2 & 2 & 7 & 0 & 0 & 0 & 0 & 0 & 0 & 0 & 10 & 1 & \dots\\
  0 & 10 & 0 & 10 & 2 & 0 & 0 & 0 & 0 & 0 & 0 & 0 & 8 & 3 & \dots\\
  0 & 2 & 8 & 9 & 3 & 0 & 0 & 0 & 0 & 0 & 0 & 0 & 7 & 5 & \dots\\
\end{tabular} }\\
Study of these and other tables led us to formulate a number of conjectures, most
of which we were able to prove. An overview of these results is
given in Theorem \ref{congthm}.
\begin{Thm} 
\label{congthm}
The following holds for the tables of $v_n$ mod $k$:\\
{\rm 1.} All $v_n$ are odd. \\
{\rm 2.} The first two rows of each table are equal.\\ 
{\rm 3.} If $k$ is even, then rows $k/2+1$ and $k/2+2$ are equal.\\
{\rm 4.} For $k$ odd, row $(k+3)/2$ contains only zeros.\\ 
{\rm 5.} For $k$ prime, the first two rows start with 1,1 followed by $k$
  occurrences of $k-1$.\\ 
{\rm 6.} For $k>2$ prime, the last $(k-1)/2$ slots of the first column vanish. \\ 
{\rm 7.} For $k>2$ prime, there is a block of zeros at the bottom (after
  $(k-1)/2$ columns), of height $(k-1)/2$ and width $(k+3)/2$.\\ 
{\rm 8.} For $k=2^q$, all rows are constant and in two fold way sweep out all
odd residues, i.e. for every odd integer $a$ with $1\le a\le 2^q$ there are precisly
two rows that have only $a$ as entry.\\    
{\rm 9.} For $k=2^q>2$ the entries in the rows $1,2,2^{q-1},2^{q-1}+1,2^{q-1}+2,2^q-1$
equal, respectively, $1,1,2^{q-1}-1,2^{q-1}+1,2^{q-1}+1,2^q-1$.\\
{\rm 10.} For $k=2^q>2$ the entries in row $a$ and row $a+2^{q-1}$ differ by
$2^{q-1}\;({\rm mod~}2^q)$.
\end{Thm}
{\it Proof}. These ten claims are proved in, respectively, lemmas \ref{lemma1},
  \ref{lemma3},  \ref{lemma3}, \ref{lemma4}, \ref{lemma6.0} \& \ref{lemma7},
  \ref{lemma8}, \ref{lemma9}, \ref{lemma13} \& \ref{equid}, \ref{nearper} and
  \ref{nearper}. \qed\\
  
\indent On computing the reductions of $v_1,\ldots,v_{32}$ modulo 32 one finds by 
part 8 of this theorem that for $k=4$ the table has constant rows $1,1,3,3$, for $k=8$ constant rows $1,1,3,3,5,5,7,7$, for $k=16$
constant rows $1,1,11,11,5,5,7,7,9,9,3,3,13,13$, $15,15$ and for $k=32$ constant rows
$1,1,27,27,21,5,7,23,9,9,19,19,29,13,31,15,17,17$, $11,11,5,21,23,7,25,25,3,3,13,29,15,31$. 
Thus, for modulus $2^q$ with $q\le 3$ we observe that pairs of values occur and that these, moreover,
are in ascending order. For $q=4$ the values still come in pairs, but the order is no longer
ascending. For $n\ge 5$ it turns out that pairs with equal values become sparser and sparser.
Notice that in the above cases for every modulus all odd values are assumed exactly
twice. By part 8 this always happens. Thus, given an odd integer $a$ and any integer
$q\ge 1$ there are infinitely many integers $m$ such that
$v_m\equiv a({\rm mod~}2^q)$ (or put more succinctly: modulo powers of two the sequence $v_n$ is
equidistributed over the odd residue classes).\\
\indent For $k$ prime, often $v_n\equiv 0({\rm mod~}k)$ for trivial reasons.  It then makes
sense to consider divisibility of $v_n$ by higher powers of $k$. Our deepest result in this
direction is provided by the following theorem.
\begin{Thm}
\label{extra1}
{\rm 1.} If $p\ge 5$ is a prime, then 
$$v_{p+3\over 2}\equiv -2p^3\;({\rm mod~}p^4){\rm ~and~}v_{p+3\over 2}\equiv 2p^3(1-p)(p-1)! 4^{p-1} 
\;({\rm mod~}p^5).$$
{\rm 2.} Let $r\ge 1$ and $p\ge 2r+1$ be a prime. 
Then
\begin{equation}
\label{starr}
v_{{p+3\over 2}+rp}\equiv C_rp^{2r+2}\;({\rm mod~}p^{2r+3}),
\end{equation}
where
$$C_r={r\over (-4)^{r-1}}\left({2r+1\over r!}\right)^2\sum_{j=0}^{2r}b_{j,r}((1-2j))_{2r-1},$$
the integers $b_{j,r}$ are defined implicitly by
$\prod_{a=1}^{2r}(2r+1-a+ax)=\sum_{j=0}^{2r}b_{j,r}x^j$,
and $((u))_a:=\prod_{j=1}^a(u+2j-2)$.

\end{Thm}
{\tt Remark}. Note that $b_{2r-j,r}=b_{j,r}$. Numerical experimentation suggests that the numerator
of $c_r$ always equals a power of 2 and that the congruence (\ref{starr}) holds for all
odd primes.\\
\indent Below we record some values of $c_r$.\\

$c_r$:
\begin{tabular}{cccc cccc}
  1 & 2 & 3 & 4 & 5 \\
  $-81$ & ${103125\over 8}$ & $-{210171535\over 64}$ & ${1308348857025\over 1024}$ & $-{11660783598520749\over 16384}$ 
\end{tabular}

\section{Proofs of the theorems}
\subsection{Some generalities}
First recall from the elementary theory of finite fields of order $p$, 
that
$$x^{p-1}-1\equiv \prod_{j=1}^{p-1} (x-j) \; ({\rm mod~}p).$$
(Here and below the letter $x$ denotes a variable.)
By substituting $x=0$ one obtains {\it Wilson's theorem}:
$$(p-1)!\equiv -1 \;({\rm mod~}p).$$
We also recall the freshman's identity $(a+b)^p\equiv a^p+b^p\;({\rm mod~}p)$, from
which we infer that if $f(x)\in \Bbb Z[x]$, then $f(x)^p\equiv f(x^p)\;({\rm mod~}p)$.
These results  will be freely used in the sequel, without further referring to them.

\begin{lem} \label{lemma0}
We have $v_n\equiv 0\;({\rm mod~}(2n-3)^2)$.
\end{lem}
{\it Proof}. The term with $j=0$ in (\ref{eq:defn}) equals $2n-3$. The
term with $j=2n-3$ equals $(2n-3)x$. 
 Hence $v_n=(2n-3)^2w_n$, where
$w_{n}=\Big[(1-x)\prod_{j=1}^{2n-4}(2n-3-j+jx)\Big]_{x^{n-2}}$. 
\qed\\

The following result was first noticed by D.~Kerner. An alternative, slightly
longer, proof was given by M.~Vlasenko.
\begin{lem}
We have $v_n\equiv 0\;({\rm mod~}(2n-3)^3)$.
\end{lem}
{\it Proof}. It suffices to show that $w_n\equiv 0\;({\rm mod~}2n-3)$. Note that, modulo 
$2n-3$, we have that
\begin{eqnarray}
w_n &\equiv & \Big[(1-x)\prod_{j=1}^{2n-4}(-j+jx)\Big]_{x^{n-2}}\equiv -(2n-4)!\Big[(x-1)^{2n-3}\Big]_{x^{n-2}}\nonumber\cr
&=& -(2n-4)!{2n-3\choose n-2}=-(2n-4)!{2n-3\over n-1}{2n-4\choose n-2}=-(2n-3)!C_{n-2},\nonumber
\end{eqnarray}
where $C_m:={1\over m+1}{2m\choose m}$ is the $m$th {\it Catalan number}. The Catalan numbers
are {\it integers} that arise in numerous counting problems.\qed\\

\noindent {\tt Remark}. It might be interesting to see whether the integers $v_n/(2n-3)^i$ 
with $i=1,2$ or 3 also have a geometric meaning.\\

The next result was obtained in collaboration with Alexander Blessing. It establisheds parts 2
and 3 of Theorem \ref{congthm}.
\begin{lem} \label{lemma3}
For $l\ge 0$ we have $v_{ln+1} \equiv v_{ln+2} \;({\rm mod~}2n)$.
\end{lem}
{\it Proof}. Since $v_1=v_2=1$ the result is trivially true for $l=0$ and thus
we may assume $l\ge 1$. We have, modulo $2n$, 
$$v_{ln+1}\equiv \Big[(1-x)\prod_{j=0}^{2ln-1}(-1-j+jx)\Big]_{x^{ln}}.$$
Furthermore, we have, modulo $2n$,
\begin{eqnarray}
v_{ln+2}&\equiv &\Big[(1-x)\prod_{j=0}^{2ln+1}(1-j+jx)\Big]_{x^{ln+1}}\equiv \Big[(1-x)x\prod_{j=0}^{2ln}
(1-j+jx)\Big]_{x^{ln+1}}\nonumber\cr
&\equiv & \Big[(1-x)\prod_{j=0}^{2ln} (1-j+jx)\Big]_{x^{ln}}
\equiv \Big[(1-x)\prod_{j=0}^{2ln} (1-(2ln-j)+(2ln-j)x)\Big]_{x^{ln}}\nonumber\cr
&\equiv & \Big[(1-x)\prod_{j=0}^{2ln}(1+j-jx)\Big]_{x^{ln}}\equiv \Big[(1-x)\prod_{j=0}^{2ln-1}(-1)(-1-j+jx)\Big]_{x^{ln}}\nonumber\cr
&\equiv & \Big[(1-x)\prod_{j=0}^{2ln-1}(-1-j+jx)\Big]_{x^{ln}}\equiv v_{ln+1}\;.
\end{eqnarray}
This concludes the proof. \qed\\

The next lemma generalizes Lemma \ref{lemma0}. It implies part 4 of Theorem \ref{congthm}.
\begin{lem} \label{lemma4}
If $k$ is odd, then $v_{lk+(k+3)/2} \equiv 0 \; ({\rm mod~}(2l+1)^2k^{2l+2})$.
\end{lem}
{\it Proof}. We have $v_{lk+(k+3)/2}\equiv [(1-x)\prod_{j=0}^{(2l+1)k}((2l+1)k-j+jx)]_{x^{lk+(k+1)/2}}$.
The terms in the product with $j=0$ and $j=(2l+1)k$ lead to a factor of
$(2l+1)^2k^2$. The remaining terms in the product that are divisible by $k$ lead to a
factor $k^{2l}$. \qed

\subsection{The sequence $\{v_n\}_{n=1}^{\infty}$ modulo primes}
The following lemma will be repeatedly used in this section.
\begin{lem} \label{lemma5}
Let $p$ be a prime and $c$ an integer. 
Then, modulo $p$,
$$ \prod_{i=1}^{p-1} (ix-i+c)
\equiv \cases{-(x-1)^{p-1} & if $p|c$,\cr
 -(x+x^2+\dots +x^{p-1})={x^p-x\over 1-x} & otherwise.}$$
\end{lem}
{\it Proof}. If $p|c$, then the result is trivial, so assume $p\nmid c$.
We can write 
$$\prod_{i=1}^{p-1} (ix-i+c)\equiv \prod_{i=1}^{p-1}i\;
\prod_{i=1}^{p-1} (x-1+c/i)\equiv {x^p-x\over 1-x},$$ where
we have used that as $i$ runs over $1,2,\cdots,p-1$, $-1+c/i$ runs over
all residues modulo $p$, except for $-1$. \qed\\

The second lemma in this section is part 1 of Theorem \ref{congthm}.
\begin{lem} \label{lemma1}
For $n\ge 1$ we have $v_n \equiv 1 \;({\rm mod~}2)$.
\end{lem} 
{\it Proof}. Modulo 2, the $2n-2$ terms in the product in
(\ref{eq:defn}) are alternatingly 1 and $x$. It thus
follows that 
$v_n \equiv [(1-x) x^{n-1}]_{x^{n-1}}\equiv 1\; ({\rm mod~}2)$.
This concludes the proof. \qed\\

The next lemma together with Lemma \ref{lemma7} establishes part 5 of Theorem \ref{congthm}.
\begin{lem} \label{lemma6.0}
Let $p$ be a prime. Then $v_{p+1}\equiv v_{p+2} \equiv 1 \;({\rm mod~}p)$.  
\end{lem}
{\it Proof}. Modulo $p$ we have
 $v_{p+1} \equiv [(1-x) \prod_{j=1}^{p-1}(-1-j-jx)^2]_{x^p} 
\equiv [(1-x) \left({x^p-x\over 1-x}\right)^2]_{x^p}$, where in the derivation of
the first congruence we noted that modulo $p$ the $j$th term
in the product (\ref{eq:defn}) is equal to the $(j+p)$th and in 
the second we used lemma~\ref{lemma5}.  
Now note that
$$
\Big[(1-x) \left({x^p-x\over
1-x}\right)^2\Big ]_{x^p}=\Big [(1-x) \left({x\over
1-x}\right)^2\Big ]_{x^p}=\Big [\sum_{k\ge 2}x^k\Big ]_{x^p}=1.
$$
Finally, by
Lemma~\ref{lemma3}, $v_{p+2}$ satisfies the same congruence as $v_{p+1}$ 
modulo $p$.
\qed

The proof of the next lemma involves congruences for binomial coefficients. 
In all cases these can be found by direct computation,
but often it is more convenient to invoke a classical result of E. Lucas. Let $n\ge m$ be
natural numbers and write $n=a_0+a_1p+a_2p^2+\cdots+a_sp^s$ and $m=b_0+b_1p+b_2p^2+\cdots+b_sp^s$ with
$0\le a_i,b_i\le p-1$. Then Lucas's theorem states that
$${n\choose m}\equiv {a_0\choose b_0}{a_1\choose b_1}\cdots {a_s\choose b_s}\;({\rm mod~}p).$$
Recall that ${a\choose b}=0$ if $b>a$.
For example, by direct computation we find that 
$${p^2-2 \choose2p-2} \equiv \left[{(-2)\cdots (-2p+1) \over 1\cdots (2p-2)}\right]'
(1-p) \equiv -(2p-1)(p-1) \equiv -1 \;({\rm mod~}p)$$ ($[\dots]'$ means skipping
multiples of $p$). By Lucas's theorem we find that
$${p^2-2 \choose2p-2}={(p-1)p+p-2\choose 1\cdot p+p-2}\equiv {p-1\choose 1}{p-2\choose p-2}\equiv
-1\; ({\rm mod~}p).$$
Likewise we immediately find using Lucas' theorem that, with $r=1$ and $p>3$,
$${2p-1\choose p-1}\equiv 1\; ({\rm mod~}p^r).$$
(This identity with $r=2$ was proved in 1819 by Charles Babbage. For $r=3$ it follows from
Wolstenholme's theorem, see F.L. Bauer \cite{Br}.)\\
\indent At various points we use the easy result that
${p-1\choose j}\equiv (-1)^j\;({\rm mod~}p)$. To see this observe that modulo $p$ the entries, 
except the two outmost ones, in the $(p+1)$th row of Pascal's triangle are zero modulo $p$. Since
each of these entries arises as sum of two elements above it in the $p$th row, the entries
in the $p$th alternate between 1 and -1. Similarly one infers that
${p-2\choose j}\equiv (-1)^j(j+1)\; ({\rm mod~}p)$.\\
\indent For a nice survey of arithmetic properties of binomial coefficients we refer the
reader to Granville \cite{AG}.
\begin{lem} \label{lemma7}
Let $p$ be a prime and $2\leq l \leq p+1$. Then $v_{lp+1} \equiv v_{lp+2}
\equiv -1 \;({\rm mod~}p)$.   
\end{lem}
{\it Proof}. We have
 $ v_{lp+1} \equiv \Big[(1-x) P_p^l(x) \Big]_{x^{lp}} \;({\rm mod~}p)$, say, with
 $P_p^l(x):=\prod_{i=1}^{p-1}(ix-i-1)^{2l}$.
As before
\begin{eqnarray}
P_p^l(x) &\equiv &(x+x^2+\dots +x^{p-1})^{2l} 
\equiv x^{2l} \Big[ (1-x^{p-1}) (1+x+x^2+\dots) \Big]^{2l} \;({\rm mod~}p)\cr
&\equiv & x^{2l} \Big[ \sum_{i=0}^{2l} x^{(p-1)i} \underbrace{ (-1)^{i} {2l
  \choose i}}_{a_{i(p-1)}}  \cdot  \sum_{n\geq 0} x^n \underbrace{2l+n
  -1 \choose n}_{b_n} \Big] \;({\rm mod~}p)\cr
&\equiv & x^{2l} \Big[\dots + x^{l(p-2)-1} \sum_{i=0}^{l-1}
a_{(l-1-i)(p-1)} b_{(i+1)(p-1)-l-1} \cr
& & +  x^{l(p-2)} \sum_{i=0}^{l-1}
a_{(l-1-i)(p-1)} b_{(i+1)(p-1)-l} +\dots \Big]\;({\rm mod~}p),\nonumber
\end{eqnarray}
where we have used that
$$\left(\sum_{k\ge 0}x^k\right)^{2l}=\sum_{n\ge 0}{2l+n-1\choose n}x^n.$$
Thus
\begin{eqnarray}
 v_{lp+1} &\equiv &  \sum_{i=0}^{l-1} a_{(l-1-i)(p-1)} (b_{(i+1)(p-1)-l} - b_{(i+1)(p-1)-l-1})\;({\rm mod~}p)\cr
& \equiv & \sum_{i=0}^{l-1} a_{(l-1-i)(p-1)} {(i+1)(p-1)+l-2 \choose 2l-2} \;({\rm mod~}p).\nonumber
\end{eqnarray}
Note that 
$$
a_{(l-1-i)(p-1)} = (-1)^{l-1-i} {2l  \choose l-1-i}=
    \cases{-2l & for $i=l-2$;\cr
1 & for $i=l-1$,}$$
and, furthermore,
$$a_{(l-1-i)(p-1)} = \cases{
~~0\;({\rm mod~}p) & for $i=0,\dots,p-l-2$;\cr 
-1\;({\rm mod~}p) & for $i=p-l-1$.}$$
The remainder of the proof requires a separate discussion for $l\leq
(p-1)/2$, for $l=(p+1)/2$, $(p+1)/2 <l <p$, for $l=p$ and for $l=p+1$.\\ 
\underline{The case $l\leq (p-1)/2$} :\\ 
Note that
\begin{eqnarray}
{(i+1)(p-1)+l-2 \choose 2l-2} &=& {[(i+1)p-i-l] \cdots
  [(i+1)p-i+l-3] \over 1 \cdots (2l-2)}\cr
  & \equiv & \cases{
0 \;({\rm mod~}p) & for $i=0,\dots,l-3$;\cr
{(-1)\cdots (-2l+2) \over  1 \cdots (2l-2)} \equiv 1 \;({\rm mod~}p) & for
$i=l-2$;\cr
{(-2)\cdots (-2l+1) \over  1 \cdots (2l-2)} \equiv 2 l-1 \;({\rm mod~}p) & for
$i=l-1$.}\nonumber
\end{eqnarray}
Hence in the sum for $ v_{lp+1}$, only the terms $i=l-2$ and $i=l-1$
contribute, that is:\\ 
$ v_{lp+1} = -2l \cdot 1 +1 \cdot (2l-1) \equiv -1 \; ({\rm mod~}p)$.\\
 \underline{The case $l=(p+1)/2$} :\\
Here we find that 
$$ {(i+1)(p-1)+l-2 \choose 2l-2} \equiv  \cases{
1 \;({\rm mod~}p) & for $i=l-2$;\cr
0 \;({\rm mod~}p) & for $0\le i\le l-1,~i\ne l-2$.}$$
Hence in the sum for $ v_{lp+1}$, only the term $i=l-2$ 
contributes and we infer that $ v_{lp+1} \equiv -2l \cdot 1\equiv -1 \;({\rm mod~}p)$.\\
\underline{The case $(p+1)/2 <l <p$} :\\
Here we compute that
$$ {(i+1)(p-1)+l-2 \choose 2l-2} \equiv  \cases{
0 \;({\rm mod~}p) & for $i=p-l,\dots,l-3$;\cr
-l \;({\rm mod~}p) & for $i=p-l-1$;\cr
2-l \;({\rm mod~}p) & for $i=l-2$;\cr
(1-2l)(l-1) \;({\rm mod~}p) & for $i=l-1$.}$$
Hence in the sum for $ v_{lp+1}$, only the terms $i=p-l-1$, $i=l-2$ and $i=l-1$
contribute. That is: $ v_{lp+1} \equiv -1 \cdot (-l) -2l \cdot (2-l) +1
\cdot (1-2l)(l-1) \equiv -1 \;({\rm mod~}p)$.\\
\underline{The case $l=p$} :\\ 
Note that $$a_{(p-1-i)(p-1)} = (-1)^{p-1-i} {2p  \choose p-1-i} = \cases{
0\;({\rm mod~}p)& for $i=0,\dots,p-2$;\cr
1\;({\rm mod~}p) & for $i=p-1$,}$$
while for $i=p-1$, $ {(i+1)(p-1)+p-2 \choose 2p-2}$ boils down to 
${p^2-2 \choose2p-2} \equiv -1 \;({\rm mod~}p)$ (see the preamble to this lemma). Hence $ v_{lp+1} \equiv
-1 \;({\rm mod~}p)$.\\
\underline{The case $l=p+1$} :\\ 
Here we find that
$$
a_{(p-i)(p-1)} = (-1)^{i+1} {2p+2 \choose p-i} \equiv \cases{
-2\;({\rm mod~}p) & for $i=0$;\cr
~~0\;({\rm mod~}p) & for $i=1,\dots, l-4$;\cr
~~1\;({\rm mod~}p) & for $i=l-3$;\cr
-2\;({\rm mod~}p) & for $i=l-2$;\cr
~~1\;({\rm mod~}p) & for $i=l-1$.}
$$
and
$$
{(i+2)(p-1) \choose 2p} \equiv \cases{
0\;({\rm mod~}p) & for $i=0$;\cr
1\;({\rm mod~}p) & for $i=l-3$;\cr
1\;({\rm mod~}p) & for $i=l-2$;\cr
0\;({\rm mod~}p) & for $i=l-1$.}
$$
Hence $ v_{lp+1} \equiv -2\cdot 0 +1\cdot 1 -2\cdot 1 +1\cdot 0 \equiv
-1 \;({\rm mod~}p)$.\\
\indent This proves that in all cases, $ v_{lp+1} \equiv -1 \;({\rm mod~}p)$. By
lemma~\ref{lemma3}, the same is true for  $ v_{lp+2}$.
\qed\\

The cases $l=p$ and $l=p+1$ in the latter proof can be proved more succinctly as is done
in the proofs of Lemma \ref{lemma6.1}, respectively Lemma \ref{lemma6.2}.
\begin{lem} \label{lemma6.1}
Let $p$ be a prime. Then $v_{p^2+1} \equiv v_{p^2+2} \equiv -1 \;({\rm mod~}p)$.  
\end{lem}
{\it Proof}. Note that, modulo $p$, the integer $v_{p^2+1}$ is congruent to 
 $$\Big [(1-x) \prod_{j=1}^{p-1}(-1-j-jx)^{2p}\Big]_{x^{p^2}}
 \equiv  \Big[(1-x) \prod_{j=1}^{p-1}(-1-j-jx^p)^2\Big]_{x^{p^2}}
 \equiv  \Big[\prod_{j=1}^{p-1}(-1-j-jy)^2\Big]_{y^p}.$$
On proceding as in the previous proof we find that
$$v_{p^2+1}\equiv \Big[\left({y^p-y\over 1-y}\right)^2\Big ]_{y^p}
=\Big [\left({y\over 1-y}\right)^2\Big ]_{y^p}
=\Big [\sum_{k\ge 1}ky^{k+1}\Big ]_{y^p}\equiv -1\;({\rm mod~}p).$$ 
Finally, by
Lemma~\ref{lemma3}, $v_{p^2+2}$ satisfies the same congruence as $v_{p^2+1}$ 
modulo $p$. \qed
\begin{lem} \label{lemma6.2}
Let $p$ be a prime. Then $v_{p^2+p+1} \equiv v_{p^2+p+2} \equiv -1 \;({\rm mod~}p)$.  
\end{lem}
{\it Proof}. We have
\begin{eqnarray}
v_{p^2+p+1}&\equiv& \Big [(1-x)(x+x^2+\cdots+x^{p-1})^{2p+2}\Big ]_{x^{p^2+p}}\; ({\rm mod~}p)\cr
&\equiv&
\Big[(1-x)(x+x^2+\cdots+x^{p-1})^2(x^p+\cdots+x^{p(p-1)})^2\Big]_{x^{p^2+p}}\; ({\rm mod~}p)\cr
&\equiv & \Big[(x-x^p)(x+x^2+\cdots+x^{p-1})\Big(\sum_{k=1}^{\infty}x^{kp}\Big)^2\Big]_{x^{p^2+p}}\;({\rm mod~}p)\cr
&\equiv & \Big[(x^2+\cdots+x^p-x^{p+1}-\cdots-x^{2p-1})\sum_{k=0}^{\infty}(k+1)x^{kp}\Big]_{x^{p^2-p}}\cr
&\equiv & \Big[\sum_{k=0}^{\infty}(k+1)x^{(k+1)p}\Big]_{x^{p^2-p}}\equiv -1\; ({\rm mod~}p).\nonumber
\end{eqnarray}
Finally, by
Lemma~\ref{lemma3}, $v_{p^2+p+2}$ satisfies the same congruence as $v_{p^2+p+1}$ 
modulo $p$. \qed\\

The next lemma establishes part 6 of Theorem \ref{congthm}.
\begin{lem} \label{lemma8}
If $p$ is an odd prime, then $v_{{p+3\over 2}+i} \equiv 0 \;({\rm mod~}p)$ for $i=0,\dots,
(p-3)/2$.
\end{lem}
{\it Proof}. In case $i=0$, the result follows by Lemma \ref{lemma0}, so assume
that $i\ge 1$. On using that modulo $p$ the $j$th term equals the $(j+p)$th term, we
find that, modulo $p$,
$$v_{{p+3\over 2}+i} \equiv \Big[4i^2(1-x) \prod_{j=1}^{p-1} (2i-j+jx) 
  \prod_{j=1}^{2i}(2i-j+jx)\Big]_{x^{(p+1)/2+j}}.$$
On invoking Lemma \ref{lemma5} and noting that $p>{p+1\over 2}+i$ we infer that
$$v_{{p+3\over 2}+i}\equiv 
\Big[4i^2(x^p-x)\prod_{j=1}^{2i}(2i-j+jx)\Big]_{x^{(p+1)/2+i}} \equiv 
\Big[-4i^2\prod_{j=1}^{2i}(2i-j+jx)\Big]_{x^{(p-1)/2+j}}.
$$
Since ${\rm deg}\bigg(\prod_{j=1}^{2i}(2i-j+jx)\bigg)=2i$ and $2i<{p-1\over 2}+i$, the 
result follows. \qed\\

The next lemma will be used in the proof of Lemma \ref{lemma9}.
\begin{lem} 
\label{slaaf}
Define $A_r(x)$ and $B_r(x)$ recursively by 
$$\cases{A_0(x)=0,~A_{r+1}(x)=(x+\ldots+x^{p-1})^r-A_r(x);\cr
B_0(x)=0,~B_{r+1}(x)=-(x+\ldots+x^{p-1})^r-B_r(x).}$$
Put
$f_r(x)=(x-1)(1+x^p+\ldots+x^{p(p-1)})(x+\ldots+x^{p-1})^{r}$.
Then 
$$f_r(x)=(-1)^r(x-1)(1+x^p+\ldots+x^{p(p-1)})+x^{p^2}A_r(x)+B_r(x),$$
where, for $r\ge 1$, the degree of $B_r(x)$ equals $(r-1)(p-1)$.
\end{lem}
{\it Proof}. Easily follows on noting that
\begin{eqnarray}
f_{r+1}(x)&=&(1+\ldots+x^{p-1})f_{r}(x)-f_r(x)\cr
&=& (x^{p^2}-1)(x+\ldots+x^{p-1})^r-f_r(x).\nonumber
\end{eqnarray}

The next lemma is part 7 of Theorem \ref{congthm}.
\begin{lem} \label{lemma9} Let $p$ be an odd prime. Suppose that  $0\le i\le (p-3)/2$
   and $(p-1)/2\le l \le p$.
   Then $v_{lp+(p+3)/2+i} \equiv 0 \; ({\rm mod~}p)$.
\end{lem}
{\it Proof}. 
Write $l=(p-1)/2+k$ (thus $k$ assumes the values $0,\ldots,(p+1)/2$).
Put $P_{2i}(x)=(2i)^{2l+2}\prod_{j=1}^{2i}(2i-j+jx)$. Note that
the degree of this polynomial equals $2i$.
Proceding as in Lemma \ref{lemma8} we infer that, modulo $p$,
\begin{eqnarray}
\label{komeropterug}
v_{lp+{p+3\over 2}+i}&\equiv & \Big[(x-1)(x+\ldots+x^{p-1})^{2l+1}P_{2i}(x)\Big]_{x^{lp+(p+1)/2+i}}\cr
&\equiv & \Big[(x-1)(x^p+x^{2p}\ldots+x^{p(p-1)})(x+\ldots+x^{p-1})^{2k}P_{2i}(x)\Big]_{x^{lp+(p+1)/2+i}}.
\end{eqnarray}
We consider the case $l=p$ (that is $k=(p+1)/2$) first. Then
\begin{eqnarray}
v_{lp+{p+3\over 2}+i}&\equiv & \Big[(x-1)(x^p+x^{2p}+\ldots+x^{p(p-1)})^2(x+\ldots+x^{p-1})P_{2i}(x)\Big]_{x^{lp+(p+1)/2+i}}\cr
&=& \Big[(x^p-x)(x^p+x^{2p}+\ldots+x^{p(p-1)})^2P_{2i}(x)\Big]_{x^{lp+(p+1)/2+i}}=0,\nonumber
\end{eqnarray}
where we used the observation that the polynomial in brackets has the form 
$$\sum_k P_{2i}(x)(c_{kp}x^{kp}+ 
c_{kp+1}x^{kp+1}),$$ and that $1+2i< (p+1)/2+i\le p-1$.\\
\indent Thus we may assume that $l\le p-1$. Notice that $lp+(p+1)/2+i<p^2$. 
Furthermore, by Lemma \ref{slaaf} we have $[B_{2k}(x)P_{2i}(x)]_{x^{(l-1)p+(p+1)/2+i}}=0$, as
$$2k-1+2i<(l-1)p+{p+1\over 2}+i.$$
Using this, (\ref{komeropterug}) and Lemma \ref{slaaf} we infer that
\begin{eqnarray}
v_{lp+{p+3\over 2}+i}&\equiv & \Big[(x-1)(x^p+\ldots+x^{p(p-1)})(x+\ldots+x^{p-1})^{2k}P_{2i}(x)\Big]_{x^{lp+(p+1)/2+i}}\cr
&= & \Big[(x-1)(x^p+\ldots+x^{p^2})(x+\ldots+x^{p-1})^{2k}P_{2i}(x)\Big]_{x^{lp+(p+1)/2+i}}\cr
&= & \Big[(x-1)(1+\ldots+x^{p(p-1)})(x+\ldots+x^{p-1})^{2k}P_{2i}(x)\Big]_{x^{(l-1)p+(p+1)/2+i}}\cr
&= & \Big[(x-1)(1+\ldots+x^{p(p-1)})P_{2i}(x)+B_{2k}(x)P_{2i}(x)\Big]_{x^{(l-1)p+(p+1)/2+i}}\cr
&= & \Big[(x-1)(1+\ldots+x^{p(p-1)})P_{2i}(x)\Big]_{x^{(l-1)p+(p+1)/2+i}}=0.\nonumber
\end{eqnarray}
This concludes the proof. \qed\\

A further question concerning the distribution of $v_n$ modulo primes is how frequently certain
residues appear. For example, is it true that the zero have density 1 ? Is it true that
the non-zero entries are equidistributed ? Questions like this can be answered for the
middle binomial coefficient ${2k\choose k}$, see e.g. \cite{BH, M1, M2}. The following lemma suggests
that perhaps techniques from the latter papers can be used to investigate this issue.
\begin{lem}
We have $v_{1+3k}\equiv v_{2+3k}\equiv {1\over k+1}{2k\choose k}\;({\rm mod~}3)$ and $v_{3k}\equiv 0\;({\rm mod~}3)$.
\end{lem}
{\it Proof}. By Lemma \ref{lemma0} we have $v_{3k}\equiv 0\;({\rm mod~}3)$. By Lemma \ref{lemma3} we
have $v_{1+3k}\equiv v_{2+3k}\;({\rm mod~}3)$.
We have, modulo 3,  
\begin{eqnarray}
v_{2+3k} &\equiv & \Big[(1-x)\prod_{j=0}^{1+6k}(1-j+jx)\Big]_{x^{1+3k}}=\Big[(1-x)(-x(1+x))^{2k} x\Big]_{x^{1+3k}}\cr
&=&\Big[(1-x)(1+x)^{2k}\Big]_{x^k}={2k\choose k}-{2k\choose k-1}={1\over k+1}{2k\choose k}.\nonumber
\end{eqnarray}
This concludes the proof. \qed

\subsection{The sequence $\{v_n\}_{n=1}^{\infty}$ modulo prime powers}
\indent The proof of the next lemma was kindly communicated to us by Carl Pomerance.
\begin{lem} 
\label{carl}
The polynomial
  $\prod_{i=0}^{p^l-1} (ix-i+j) \;({\rm mod~}p^l)$, as a polynomial in $x$,
  depends only on the class $j 
  \;({\rm mod~}p)$ (ie. replacing $j$ by $j+kp$ would yield the same result).  
\end{lem}
{\it Proof}. Let $f_j(x) = \prod_{i=0}^{p^{r}-1}(ix+j)$.
If $p|j$, then there are $p^{r-1}$ factors divisible by $p$
and $p^{r-1} \ge r$, so that $f_j(x)\equiv 0 \;({\rm mod~}p^{r})$.
So assume $p\nmid j$.  Let $k$ be the inverse of $j$,
so $jk \equiv 1\;({\rm mod~}p^{r})$.  Then modulo $p^{r}$, we have
$f_j(x) \equiv j^{p^r} f_1(x)$  (since the expression $ik$ runs
over a complete residue system modulo $p^{r}$ as $i$ runs).
Now say $j \equiv j_1 \;({\rm mod~}p)$, say $j_1= j+kp$. Using
induction with respect to $r$ one then easily sees that
$j_1^{p^r}=(j+kp)^{p^r}\equiv j^{p^r}\;({\rm mod~}p^r)$, and we
are done.\qed\\

\noindent {\it Proof of Theorem} \ref{extra1}. Part 1. We have
\begin{eqnarray}
v_{p+3\over 2}&=& p^2\Big[(1-x)\prod_{j=1}^{p-1}(p-j+jx)\Big]_{x^{p-1\over 2}}\cr
&\equiv & p^2\Big[(1-x)\prod_{j=1}^{p-1}(-j+jx)+(1-x)p\sum_{k=1}^{p-1}\prod_{j=1\atop j\ne k}^{p-1}
(-j+jx)\cr
&~& + (1-x)p^2\sum_{1\le k<r\le p-1}\prod_{j=1\atop j\ne k,r}^{p-1}(-j+jx)\Big]_{x^{p-1\over 2}}\;({\rm mod~}p^5)\cr
&\equiv & -p^2(p-1)!\Big[(x-1)^p+(x-1)^{p-1}p\sum_{k=1}^{p-1}{1\over k}\cr
&~&+(x-1)^{p-2}p^2\sum_{1\le k<r\le p-1}{1\over kr}\Big]_{x^{p-1\over 2}}\;({\rm mod~}p^5)\cr
&\equiv & -p^2\{(p-1)!\}\Big[(x-1)^p\Big]_{x^{p-1\over 2}}\;({\rm mod~}p^5),\cr
&\equiv & -p^2\{(p-1)!\}{p\choose {p-1\over 2}}(-1)^{p+1\over 2}\;({\rm mod~}p^5),\cr
&\equiv & 2p^3\{(p-1)!\}(1-p){p-1\choose {p-1\over 2}}(-1)^{p-1\over 2}\;({\rm mod~}p^5),\nonumber
\end{eqnarray}
where we used that $\sum_{k=1}^{p-1}1/k\equiv 0\;({\rm mod~}p^2)$ (this is Wolstenholme's theorem \cite[Theorem 115]{HW})
and $\sum_{1\le k<r\le p-1}{1/kr}\equiv 0\;({\rm mod~}p)$. To see the latter congruence note that
$$(p-1)!\sum_{1\le k<r\le p-1}{1\over kr}=\Big[\prod_{j=1}^{p-1}(x-j)\Big]_{x^{p-3}}\equiv \Big[x^{p-1}-1\Big]_{x^{p-3}}=0\;({\rm mod~}p).$$
Now it is an easy consequence of Eisenstein's congruence (1859), see \cite[Theorem 132]{HW}, which 
states that
$${2^{p-1}-1\over p}\equiv 1+{1\over 3}+{1\over 5}+\cdots+{1\over p-2}\;({\rm mod~}p),$$
that (see ibid. Theorem 133)
${p\choose (p-1)/2}(-1)^{p-1\over 2}\equiv 4^{p-1}\;({\rm mod~}p^2)$.
(Indeed, by Morley's congruence (1895), cf. \cite{Cai},  this congruence is even valid modulo $p^3$.)
We thus finally infer that  $v_{p+3\over 2}\equiv 2p^3(1-p)\{(p-1)!\} 4^{p-1} \;({\rm mod~}p^5)$, which
of course implies that $v_{p+3\over 2}\equiv -2p^3\;({\rm mod~}p^4)$.\\
\indent Part 2. We have the formal series identity
$${1\over (1-x)^{2r}}=\sum_{k=0}^{\infty}{k+2r-1\choose 2r-1}x^k.$$
Note that
$$[(x-1)^{-2r}]_{x^{{p-1\over 2}-s+rp}}\equiv {(1-2s)_{2r-1}\over 2^{2r-1}(2r-1)!}\;({\rm mod~}p).$$
Using the latter congruence we find that, modulo $p^{2r+3}$.
\begin{eqnarray}
v_{{p+3\over 2}+rp}&=&\Big[(1-x)\prod_{j=0}^{(2r+1)p}((2r+1)p-j+jx)\Big]_{x^{{p+1\over 2}+rp}}\cr
&\equiv &p^{2r+2}\Big[(1-x)\prod_{j=0}^{2r+1}(2r+1-j+jx)\prod_{j=0\atop p\nmid j}^{(2r+1)p}(-j+jx)\Big]_{x^{{p+1\over 2}+rp}}
\cr
&\equiv &(2r+1)^2p^{2r+2}\Big[\Big(\sum_{j=0}^{2r}b_{j,r}x^j\Big)(x-1)^{(2r+1)p-2r}\Big]_{x^{(p-1)/2+rp}}\cr
&\equiv &(2r+1)^2p^{2r+2}\Big[\Big(\sum_{j=0}^{2r}b_{j,r}x^j\Big)(x-1)^{(2r+1)p}\sum_{k=0}^{\infty}{k+2r-1\choose 2r-1}x^k\Big]_{x^{(p-1)/2+rp}}\cr
&\equiv &-(2r+1)^2p^{2r+2}\Big[\Big(\sum_{j=0}^{2r}b_{j,r}x^j\Big)\sum_{l=0}^rx^{lp}(-1)^l\sum_{k=0}^{\infty}{k+2r-1\choose 2r-1}x^k\Big]_{x^{(p-1)/2+rp}}\cr
&\equiv & -(2r+1)^2p^{2r+2}\sum_{l=0}^r{2r+1\choose l}(-1)^l\sum_{j=0}^{2r}b_{j,r}{{p-1\over 2}-j+(r-l)p+2r-1\choose 2r-1}\cr
&\equiv & -{(2r+1)^2p^{2r+2}\over 2^{2r-1}(2r-1)!}\sum_{l=0}^r{2r+1\choose l}(-1)^l\sum_{j=0}^{2r}b_{j,r}((1-2j))_{2r-1}\cr
&\equiv & {(-1)^{r-1}(2r+1)^2\over 2^{2r-1}(2r-1)!}{2r\choose r}p^{2r+2}\sum_{j=0}^{2r}b_{j,r}((1-2j))_{2r-1}\cr
&=& C_r p^{2r+2},\nonumber
\end{eqnarray}
where in the one but last step we used the identity
$$(-1)^r{2r\choose r}=\sum_{l=0}^r{2r+1\choose l}(-1)^l,$$
which is obtained by comparing the coefficient of $x^r$ of both sides
of the identity\\ $(1-x)^{-1}(1-x)^{2r+1}=(1-x)^{2r}$.
This finishes the proof. \qed

\subsection{The sequence $\{v_n\}_{n=1}^{\infty}$ modulo powers of two}
Before we can consider the sequence modulo powers of two we need some
preparatory lemmas.
\begin{lem} \label{lemma12}
If $j$ is odd, then $\prod_{i=0}^{2^q-1} (ix-i+j)^2 \equiv x^{2^q} \;({\rm mod~}2^q)$.
\end{lem} 
{\it Proof}. By induction with respect to $q$. For $q=1$ the result 
is obvious. Assume the result is established for $1\le q\le q_1$.
We write $$\prod_{i=0}^{2^{q_1+1}-1} (ix-i+j)^2=
\prod_{i=0}^{2^{q_1}-1} (ix-i+j)^2  \prod_{i=2^{q_1}}^{2^{q_1+1}-1} (ix-i+j)^2
=P_1(x)P_2(x),$$ say.
Note that $P_1(x)\equiv P_2(x)\;({\rm mod~}2^{q_1})$. The induction hypothesis
thus implies that we can write $P_1(x)=x^{2^{q_1}}+2^{q_1}f_1(x)$ and  
$P_2(x)=x^{2^{q_1}}+2^{q_1}f_2(x)$. Since 
$(ix-i+j)^2\equiv ((i+2^{q_1})x-(i+2^{q_1})+j)^2 \;({\rm mod~}2^{q_1+1})$, it even
follows that $P_1(x)\equiv P_2(x)\;({\rm mod~}2^{q_1+1})$, from which we infer
that $f_1(x)\equiv f_2(x)\;({\rm mod~}2)$ and
hence $f_1(x)+f_2(x)\equiv 0\;({\rm mod~}2)$. It follows that modulo $2^{q_1+1}$ the
product under consideration equals
$$P_1(x)P_2(x)=(x^{2^{q_1}}+2^{q_1}f_1(x))(x^{2^{q_1}}+2^{q_1}f_2(x))
=x^{2^{q_1+1}} \;({\rm mod~}2^{q_1+1}).$$
This concludes the proof. \qed\\

In the course of  the above proof we have showed that
$$\prod_{i=0}^{2^q-1} (ix-i+j)^2  \equiv \prod_{i=2^q}^{2^{q+1}-1} (ix-i+j)^2 \;\;({\rm mod~}2^{q+1}).$$
The next result shows that the same identity holds true for the `square roots'. Using this
the `square root' of the left hand side of Lemma \ref{lemma12} can be computed (Lemma \ref{lemma13b}).
\begin{lem}
\label{lemma13a} Let $j$ be odd and $q\ge 2$. Then
$$\prod_{i=0}^{2^q-1} (ix-i+j) \equiv \prod_{i=2^q}^{2^{q+1}-1} (ix-i+j)\;\;({\rm mod~}2^{q+1}).$$
\end{lem}
{\it Proof}. 
It is an easy observation that, modulo 2,
we have for $0\le k\le 2^q-1$ that
$$\prod_{a=0,~a\ne k}^{2^q-1}(j-a+ax)\equiv \cases{x^{2^{q-1}-1} & if $k$ is odd;\cr
x^{2^{q-1}} & if $k$ is even.}$$
Using this identity we find that, modulo $2^{q+1}$, 
\begin{eqnarray}
\prod_{i=2^q}^{2^{q+1}-1} (ix-i+j)&=&\prod_{i=0}^{2^{q}-1} (ix-i+j+2^q(x-1))\cr
&\equiv &\prod_{i=0}^{2^q-1} (ix-i+j)+2^q(x-1)\sum_{k=0}^{2^q-1}\prod_{i=0,~i\ne k}^{2^q-1} (ix-i+j)\cr
&\equiv &\prod_{i=0}^{2^q-1} (ix-i+j)+2^q(x-1)\Big(x^{2^{q-1}}\sum_{2|k}^{2^q-2}1
+x^{2^{q-1}-1}\sum_{2\nmid k}^{2^q-1}1\Big)\cr
&\equiv &\prod_{i=0}^{2^q-1} (ix-i+j)+2^q(x-1)(x^{2^{q-1}}2^{q-1}
+x^{2^{q-1}-1}2^{q-1})\cr
&\equiv &\prod_{i=0}^{2^q-1} (ix-i+j).\nonumber
\end{eqnarray}
This finishes the proof. \qed

\begin{lem} \label{lemma13b} Let $j$ be odd and $q\ge 3$. We have
$$\prod_{i=0}^{2^q-1} (ix-i+j) \equiv x^{2^{q-1}-2}\Big[2^{q-1}(x^4+x^3+x+1)+x^2\Big] 
\; \;({\rm mod~}2^{q}).$$
\end{lem}
{\it Proof}. Similar to that of Lemma \ref{lemma12}, but with the difference
that instead of the equality $P_1(x)\equiv P_2(x)\;({\rm mod~}2^{q+1})$, we
use Lemma \ref{lemma13a}. \qed\\

\noindent {\tt Remark}. By Lemma \ref{carl} it suffices to work in the proofs of Lemma \ref{lemma12}, \ref{lemma13a}
and \ref{lemma13b} with $j=1$. \\

The next result established part of parts 8 and 9 of Theorem \ref{congthm}.
\begin{lem} \label{lemma12a}
For $q\ge 1$ we have $v_{2^q}\equiv -1\;({\rm mod~}2^q)$.
\end{lem}
{\it Proof}. Put $P_q(x)=(1-x)\prod_{j=0}^{2^{q+1}-3}(-3-j+jx)$. We want to
compute the coefficient of $x^{2^q-1}$ in $P_q(x)$ modulo $2^q$. On invoking
Lemma \ref{lemma12} one finds that
$$P_q(x)(1+2x)(2+x)\equiv x^{2^q}(1-x)\;({\rm mod~}2^q),$$
from which we infer that
\begin{eqnarray}
P_q(x) & \equiv & x^{2^q-q}(1-x)\sum_{k=0}^{\infty}(-2)^kx^k\sum_{r=0}^{q-1}(-2)^{q-1-r}x^r\;({\rm mod~}2^q)\nonumber\\
&\equiv &x^{2^q-q}(1-x)\sum_{m=0}^{2q-2}a_mx^m \equiv 
x^{2^q-q}\sum_{m=0}^{2q-1}b_mx^m\;({\rm mod~}2^q),\nonumber
\end{eqnarray}
where
$$a_m\equiv  \cases{-(-2)^{q-1-m}/3 & if $0\le m\le q-1$;\cr -(-2)^{-q+1+m}/3 & if $q\le m\le 2q-2$}
{~and~}b_m\equiv \cases{-(-2)^{q-1-m} & if $0\le m\le q-1$;\cr (-2)^{m-q} & if $q\le m\le 2q-1$.}$$
Thus $v_{2^q}\equiv b_{q-1}\equiv -1\;({\rm mod~}2^q)$. \qed\\

Recall that we defined $v_0=-1$. The reason for this is that this definition allows us
to also formulate the next lemma, which together with Lemma \ref{equid} is part 8 
of Theorem \ref{congthm}, with $j=0$.
\begin{lem} \label{lemma13} (Periodicity.)
Suppose that $i,k\ge 0$. We have $v_{k 2^q+i}\equiv v_i \;({\rm mod~}2^q)$.
\end{lem}
{\it Proof}. First assume that $i\ge 2$. Note that
  $$v_{k 2^q+i} \equiv \Big[ (1-x) \prod_{j=0}^{2^q-1}
  (2i-3-j+jx)^{2k} \prod_{j=0}^{2i-3} (2i-3-j+jx) \Big]_{x^{k2^q+i-1}}
  \;({\rm mod~}2^q).$$  By lemma~\ref{lemma12}, the first product
  equals $x^{k2^q}$ mod $2^q$.  Thus  
  $$ v_{k 2^q+i} \equiv \Big[(1-x) \prod_{j=0}^{2i-3} (2i-3-j+jx)\Big]_{x^{i-1}}\equiv v_i \;({\rm mod~}2^q).$$
In order to deal with the case $i=1$, we note that, 
using Lemma \ref{lemma3}, $v_{k2^q+1}\equiv v_{k2^q+2}\equiv v_2\equiv v_1\; ({\rm mod~}2^q)$. 
In case $i=0$ one finds proceding as above that, for $k\ge 1$, $v_{k2^q}\equiv v_{2^q}\; ({\rm mod~}2^q)$.
On invoking Lemma \ref{lemma12a} it then follows that $v_{k2^q}\equiv v_{2^q}\equiv v_0\; ({\rm mod~}2^q)$. 
\qed\\

The next result yields a part of part 9 of Theorem \ref{congthm}.
\begin{lem}
\label{plusnul} 
Suppose that $q\ge 1$. Then $v_{2^{q-1}}\equiv 2^{q-1}-1 \;({\rm mod~}2^q)$.
\end{lem}
{\it Proof}. Similar to that of Lemma \ref{lemma12a}. For $q\le 3$ one verifies the
claim numerically. So assume $q\ge 4$. We want to compute the coefficient of
$x^{2^{q-1}-1}$ in $P_{q-1}(x)$ modulo $2^q$. On invoking Lemma \ref{lemma13b} one
finds that
$$P_{q-1}(x)(1+2x)(2+x)\equiv x^{2^{q-1}-2}\Big[2^{q-1}(x^4+x^3+x+1)+x^2\Big]\;({\rm mod~}2^q),$$
whence
$P_{q-1}(x)\equiv x^{2^{q-1}-q-2}\Big(\sum_{m=0}^{2q-1}b_mx^m\Big)
\Big[2^{q-1}(x^4+x^3+x+1)+x^2\Big]\;({\rm mod~}2^q)$.
(Note that the assumption $q\ge 4$ implies that 
$2^{q-1}-q-2\ge 0$.) Thus modulo $2^q$ the coefficient of $x^{2^{q-1}-1}$, that
is $v_{2^{q-1}}$, equals
\begin{eqnarray}
v_{2^{q-1}}&\equiv & b_{q-1}+2^{q-1}(b_{q-3}+b_{q-2}+b_q+b_{q+1})\cr
&\equiv & -1+2^{q-1}(-4+2-2+1)\equiv 2^{q-1}-1.\nonumber
\end{eqnarray}
This completes the proof. \qed \\

The next result with $i=0,1$ and 2 yields a part of part 9 of Theorem \ref{congthm}. It also
yields part 10 of Theorem \ref{congthm}.
\begin{lem}
\label{nearper}
For $i\ge 0$ and $q\ge 2$ we have $v_{2^{q-1}+i}\equiv v_i+2^{q-1} \; ({\rm mod~}2^q)$.
\end{lem}
{\it Proof}. For $q=2$ one checks the result numerically and so we may assume $q\ge 3$.
For $i=0$ the result follows by Lemma \ref{plusnul}. 
Note that, a priori, $v_{2^{q-1}+i}\equiv v_i\; ({\rm mod~}2^{q-1})$ and so
either $v_{2^{q-1}+i}\equiv v_i\;({\rm mod~}2^{q})$ or $v_{2^{q-1}+i}\equiv v_i+2^{q-1}\; ({\rm mod~}2^{q})$.  
Let us first assume that $i\ge 2$. The idea of the proof is to use Lemma \ref{lemma13a} to write
$v_{2^{q-1}+i}\equiv v_i+2^{q-1}[f_{q,i}(x)]_{x^{2^{q-1}+i-1}} \; ({\rm mod~}2^q)$. An easy
computation then shows that $[f_{q,i}(x)]_{x^{2^{q-1}+i-1}}$ is odd, thus finishing the proof.\\
\indent More precisely, one first notes that, modulo $2^q$, 
\begin{eqnarray} 
v_{2^{q-1}+i}& \equiv &\Big[(1-x)\prod_{j=0}^{2^q-1}(2i-3-j+jx)\prod_{j=2^q}^{2^q+2i-3}(2i-3-j+jx)\Big]_{x^{2^{q-1}+i-1}}\cr
&\equiv &\Big[(1-x)\prod_{j=0}^{2^q-1}(2i-3-j+jx)\prod_{j=0}^{2i-3}(2i-3-j+jx)\Big]_{x^{2^{q-1}+i-1}}\cr
&\equiv & \Big[(1-x)x^{2^{q-1}-2}(2^{q-1}(x^4+x^3+x+1)+x^2)\prod_{j=0}^{2i-3}(2i-3-j+jx)\Big]_{x^{2^{q-1}+i-1}}\cr
&\equiv & v_i + 2^{q-1}\Big[(1-x)(x^4+x^3+x+1)\prod_{j=0}^{2i-3}(2i-3-j+jx)\Big]_{x^{i+1}}\cr
&\equiv & v_i + 2^{q-1}\Big[(1-x)(x^4+x^3+x+1)x^{i-1}\Big]_{x^{i+1}}\nonumber
\end{eqnarray}
Since $[(1-x)(x^4+x^3+x+1)x^{i-1}]_{x^{i+1}}=[(1-x)(x^4+x^3+x+1)]_{x^2}=-1$ is odd, the
result follows for $j\ge 2$. On combining this result for $i=2$ and Lemma \ref{lemma3} we find that
$v_{2^{q-1}+1}\equiv v_{2^{q-1}+2}\equiv v_2=v_1\; ({\rm mod~}2^q)$. Thus the result follows for
every $j\ge 0$. \qed\\

Using induction and the latter lemma one then easily infers the following result which, together
with Lemma \ref{lemma13}, gives part 8 of Theorem \ref{congthm}.
\begin{lem} 
\label{equid}
(Equidistribution.)
Let $q\ge 1$. For every odd integer $a$ there are precisely two integers $1\le j_1<j_2\le 2^q$
such that $v_{j_1}\equiv a \; ({\rm mod~}2^q)$ and $v_{j_2}\equiv a\; ({\rm mod~}2^q)$.
\end{lem}

\subsection{On a result of Paolo Dominici}
Let $S_k(x_1,\ldots,x_r)$ denote the $k$th elementary symmetric function
in $r$ variables, i.e. $S_1(x_1,\ldots,x_r)=x_1+\ldots+x_r$, 
$S_2(x_1,\ldots,x_r)=x_1x_2+x_1x_3+\ldots+x_{r-1}x_r$, etc.. 
Paolo Dominici \cite{Dom} states the following result for $v_n$ without
reference.
\begin{Thm}
\label{paolo}
For $1\le i\le 2n-4$ we put $y_i=i/(2n-3-i)$.
Then
$$v_n=(2n-3)^2 (2n-4)!\{S_{n-2}(y_1,\ldots,y_{2n-4})-S_{n-1}(y_1,\ldots,y_{2n-4})\}.$$
\end{Thm}
We will now derive this result from (\ref{eq:defn}). We need two lemmas
\begin{lem}
\label{lemm1}
Let $L_1(x),\ldots,L_r(x)$ be linear polynomials, then
$${1\over m!}{d^m\over dx^m}\{L_1(x)\cdots L_r(x)\}=S_m\left({L_1'(x)\over L_1(x)},
\ldots,{L_r'(x)\over L_r(x)}\right)L_1(x)\cdots L_r(x).$$
\end{lem}
Another observation we need is the following:
\begin{lem} 
Let $x_1,\ldots,x_r$ be distinct non-zero elements such
that $x_1\cdot x_2\cdots x_r=1$ and $\{x_1,\ldots,x_r\}=\{{1\over x_1},\ldots,{1\over x_r}\}$. Then
$S_{r-k}(x_1,\ldots,x_r)=S_k(x_1,\ldots,x_r)$, with $1\le r\le k$.
\end{lem}
{\it Proof}. Note that 
$S_{r-k}(x_1,\ldots,x_r)=S_k({1\over x_1},\ldots,{1\over x_r})x_1\cdots x_r=S_k(x_1,\ldots,x_r)$,
where in the derivation of the first equality we used the assumption that $x_i\ne 0$ and
in that of the second the remaining assumptions. \qed
\begin{cor}
\label{siem}
For $1\le k\le 2n-5$ we have $S_k(y_1,\ldots,y_{2n-4})=S_{2n-4-k}(y_1,\ldots,y_{2n-4})$.
\end{cor}

\noindent {\it Proof of Theorem} \ref{paolo}. Put $P_n(x)=\prod_{j=1}^{2n-4}(2n-3-j+jx)$.
By definition we have
$$v_n=\Big[(1-x)\prod_{j=0}^{2n-3}(2n-3-j+jx)\Big]_{x^{n-1}}=(2n-3)^2\Big[(1-x)P_n(x)\Big]_{x^{n-2}}.$$
Thus, 
\begin{equation}
\label{v1}
v_n=(2n-3)^2\{ [P_n(x)]_{x^{n-2}}-[P_n(x)]_{x^{n-3}}\}.
\end{equation}
On noting that $[P_n(x)]_{x^m}={1\over m!}{d^m\over dx^m}P_n(x)\big|_{x=0},$
we obtain on invoking Lemma \ref{lemm1} that
\begin{equation}
\label{gadoor}
[P_n(x)]_{x^m}=(2n-4)!S_m(y_1,\ldots,y_i,\ldots,y_{2n-4}).
\end{equation}
Combining (\ref{gadoor}) with (\ref{v1}) yields that 
$$v_n=(2n-3)^2(2n-4)!\big\{S_{n-2}(y_1,\ldots,y_{2n-4})
-S_{n-3}(y_1,\ldots,y_{2n-4})\big\},$$
or, on invoking Corollary \ref{siem},
$$v_n=(2n-3)^2(2n-4)!\big\{S_{n-2}(y_1,\ldots,y_{2n-4})
-S_{n-1}(y_1,\ldots,y_{2n-4})\big\},$$
This concludes the proof. \qed\\

\indent {}From the above proof we infer that we may alternatively define $v_n$ by
\begin{equation}
\label{alterego}
v_n=\Big[(x-1)\prod_{j=0}^{2n-3}(2n-3-j+jx)\Big]_{x^{n}}.
\end{equation}
We leave it to the reader to use the observation that
$P(x):=\prod_{j=1}^{2n-4}(2n-3-j+jx)$ is selfreciprocal, i.e. satisfies
$P(1/x)x^{2m-4}=P(x)$ to infer (\ref{alterego}) directly from (\ref{eq:defn}).\\ 

\indent Many of the congruences can be also proved using Theorem \ref{paolo}. As an example
we will show that if $p$ is an odd prime, then $v_{3(p+1)/2}\equiv -81p^4\;({\rm mod~}p^5)$.
This is the case $r=1$ of part 2 of Theorem \ref{extra1}.\\

\noindent {\it Proof of part 2 of Theorem \ref{extra1} in case 
$r=1$.} Set $n=3(p+1)/2$. Note that $(2n-3)^2(2n-4)!\equiv -18p^4\;({\rm mod~}p^5)$. It thus remains
to be proven that the expression in braces in Theorem \ref{paolo} equals $9/2$ modulo $p$.
It turns out to be a little easier to work with $w_i=-y_i$. Note that
$(-1)^rS_r(w_1,\ldots,w_{3p-1})=S_r(y_1,\ldots,y_{3p-1})$. We have
$w_i=i/(i-3p)$ for $1\le i\le 3p-1$. Thus $w_p=-1/2$, $w_{2p}=-2$ and the remaining
$w_i$ satisfy $w_i\equiv 1\;({\rm mod~}p)$. Hence
$S_r(w_1,\ldots,w_{3p-1})\equiv S_r(-1/2,-2,1,1,\ldots,1)\;({\rm mod~}p)$, where $2\le r\le 3p-1$.
In the symmetric function $S_r(z_1,\ldots,z_{3p-1})$ there are ${3p-3\choose r}$ terms
containing neither $z_1$ nor $z_2$. 
There are ${3p-3\choose r-1}$ terms containing $z_1$, but not $z_2$. Finally there
are ${3p-3\choose r-2}$ terms containing both $z_1$ and $z_2$.
It follows that,  modulo $p$,
$$(-1)^rS_r(y_1,\ldots,y_{3p-1})\equiv S_r(-{1\over 2},-2,1,\ldots,1)={3p-3\choose r}-(2+{1\over 2}){3p-3\choose r-1}+
{3p-3\choose r-2}.$$
Modulo $p$ we have
\begin{eqnarray}
&(-1)^n&\{S_{n-2}(y_1,\ldots,y_{2n-4})- S_{n-1}(y_1,\ldots,y_{2n-4})\}\cr
&\equiv &{3p-2\choose n-1}-{5\over 2}{3p-2\choose n-2}+{3p-2\choose n-3}\cr
&\equiv &2{p-2\choose n-p-1}-{5}{p-2\choose n-p-2}+2{p-2\choose n-p-3}\cr
&\equiv &\left[2{p-2\choose n-p-1}+4{p-2\choose n-p-2}+2{p-2\choose n-p-3}\right]-{9}{p-2\choose n-p-2}\cr
&\equiv &\left[2{p\choose n-p-1}\right]-{9}{p-2\choose n-p-2}\cr
&\equiv &-{9}{p-2\choose n-p-2}\equiv (-1)^{n}{9}(n-p-1)\equiv (-1)^n{9\over 2}.\nonumber
\end{eqnarray}
This completes the proof. \qed
 
\section{Asymptotics}

Given a sequence of coefficients, there are many things we would like
to know about it.  Apart from the search for a generating function and
for a recursion formula, an interesting question is the asymptotic
behaviour.  We remind the reader that candidate Fourier series for
modular forms of weight $2k$ for $SL(2,\Z)$ must have coefficients
growing like $n^{2k-1}$ (and $n^k$ for cusp forms).

In our case, without prior knowledge of the alternative definition
(\ref{eq:defn}), we only managed to compute the first 80 values of
$v_n$ using the {\em Schubert} package for intersection theory
\cite{KS-92} or the first 225 values using the rational function (with
the dummy variables $w_i$).  Numerically, it is readily seen that the
leading term for the $v_n$ is $e^{2n\log n}$.  As this is strongly
reminiscent of the behaviour of $(2n)! = \exp(2n\log(2n) -2n
  +\half\log 2n +\half\log 2\pi +O({1\over n}))$, we rather study the
behaviour of $\log \frac{v_n}{(2n)!}$ and find now the leading term to
be $2n$.  Subtracting it, we find the next-to-leading term to be $-4
\log n$, easily verified by applying $n\p_n$ (ie. taking subsequent
differences and multiplying by $n$).  The next term is a constant,
$C=-5.62...$, which we find difficult to recognize.  We have learnt
from Don Zagier a smart technique which enables to determine a large
number of digits of $C$; we present it below under the name of
asymp$_k$ trick.

\subsection{The asymp$_k$ trick}

Assume we are given numerically a few hundred terms of a sequence $\,s=\{s_n\}_{n\in\N}$
which we believe has an asymptotic expansion goes in inverse powers of $n$, ie. 
$$s_n \;\sim\; c_0 + {c_1\over n}  +{c_2\over n^2}  +\dots .$$ \\
{\underline{Goal}}: determine the coefficients $c_i$ numerically.\\
{\underline{Trick}}: Choose some moderate value of $k$ (say $k=8$) and define a new sequence
$s^{(k)}$ as $\displaystyle{1\over k!}\Delta^kN^ks$, where $\Delta$ is the difference operator
$(\Delta u)_n=u_n-u_{n-1}$ and $N$ the multiplication operator $(Nu)_n=nu_n$, i.e., 
$$s_n^{(k)}\;=\; \sum_{j=0}^k\frac{(-1)^j}{j!\,(k-j)!}\,(n-j)^k\,s_{n-j}\;.$$
For $n$ large we have (assuming the above asymptotic expansion for $s$ itself)
$$s_n^{(k)}\;=\;c_0 \,+\,(-1)^k\,{c_{k+1}\over n^{k+1}}
  + \,(-1)^k\,{\bigl((k+1)c_{k+2}-{k+1\choose 2}c_{k+1}\bigr)\over n^{k+2}} \,+\,\dots\,.$$
Thus, while $s_n$ approximates $c_0$ only to within an accuracy $O(n^{-1})$, $s_n^{(k)}$ approximates
it to the much better accuracy $O(n^{-k})$, so we obtain a very good approximation for $c_0$. 
Call this operation \underline{asymp$_k$}.  The further coefficients $c_i$ are then obtained inductively:
if $c_0,\dots,c_{i-1}$ are known to high precision, we get $c_i$ by applying asymp$_k$ to the sequence 
$n^i\bigl(s_n-c_0-\cdots-c_{i-1}/n^{i-1})=c_i+c_{i+1}/n+\cdots\,$.

The crucial point in the success of asymp$_k$ is that the operator $\Delta^k$ sends $n^k$ to
$k!$ and kills polynomials of degree $<k$, so that all the intermediate terms of the expansion
of $s_n$ between $c_0$ and $c_kn^{-k}$ disappear.

Variants of asymp$_k$ allow one to deal for example with asymptotic expansions of the form
\begin{enumerate} 
\item[(I)] $s_n\,\sim\, A\log n+ c_0 + {c_1/n}  +{c_2 /n^2}  +\cdots $
\item[(II)] $s_n\,\sim\, Bn+A\log n+ c_0 + {c_1/n}  +{c_2 /n^2}  +\cdots $ 
\item[(III)] $s_n\,\sim\,An^{\lambda}(1+{c_1/n}+{c_2/n^2}+\cdots )$
\end{enumerate} 
In case (I) we can apply asymp$_k$ to the sequence $n(s_{n+1}-s_n)$, which has 
the form $A+c_1'/n+c_2'/n^2+\cdots $, to obtain $A$ to high precision, after which we 
apply the original method to $\{s_n-A\log n\}$.  In case~(II) we apply  asymp$_k$
twice to $\,\Delta s\,$ to get $B$ and~$A$, and then subtract (our approximation for) $B\log n+A$
from $s_n$ and apply the standard version.  For case~(III) we can either look at $\{\log s_n\}$ and
apply variant~(I) or else apply asymp$_k$ to $\{n(s_{n+1}/s_n-1)\}$ to get $\lambda$
and then apply the standard method to $\{s_n/n^\lambda\}$.\smallskip

\noindent {\tt Remark 1}. Applying the operation asymp$_k$ with suitably chosen $k$ gives a rapidly convergent
sequence $s^{(k)}\,$.  To estimate how many decimals are probably correct, we look at some relatively widely 
spaced elements of this sequence (e.g., the terms $s_n^{(k)}$ with $n=300$, 400, 500 if we know 500 terms
of the sequence $s$) and see how many of their digits agree.\smallskip

\noindent {\tt Remark 2}. One also has to experiment to find the optimal choice of $k$. Typically one uses
$k=5$ if one knows 200 terms of $s$ and $k=8$ if one knows 1000 terms. This suggests that perhaps 
$k\approx\log N$ is a good choice for a generic sequence with $N$ computed terms.\smallskip

\noindent {\tt Remark 3}. The asymp$_k$ trick was first described in a paper of Zagier \cite{Z}. Here he
considers the Stoimenov numbers $\xi_D$ which bound the number $V(D)$ of linearly indepedent 
Vassiliev invariants of degree $D$. Stoimenov himself thought that $\xi_D$ behaves `something
like $D!/1.5^D$'. Calculating the values up to $D=200$ and applying a variation of asymp$_k$ suggested
an asymptotic formula of the form
$$\xi_D\;\sim\;{D!\sqrt{D}\over (\pi^2/6)^D}\,\Big(C_0+{C_1\over D}+{C_2\over D^2}+\cdots \Big)\,,$$
with $C_0\approx 2.704332490062429595$, $C_1\approx -1.52707$ and $C_2\approx -0.269009$.
Subsequently Zagier was able to prove this with explicitly computable constants $C_i$. In particular,
$C_0=12\sqrt{3}\pi^{-5/2}e^{\pi^2/12}$, which agrees to the accuracy given above with the empirically
obtained value.

\subsection{Application to the asymptotics of $v_n$}

In our case of sequence $v_n$ of lines in a hypersurface of $\P^n$, the
coefficients $c_0=:C$ is difficult to recognize, but all other
coefficients, $c_1, c_2,\dots$ are rational numbers which we easily
recognize from a sufficient number of digits.  Once the first few
rational coefficients have been found and the corresponding terms
subtracted from the sequence $s$, the constant term $C$ can be
obtained with 30 digits, say.  This is enough to feed to the PARI
software and apply the function {\tt
  lindep([C,1,log(Pi),log(2),log(3)])} to find a rational linear
combination of $C$ in terms of a given basis (educated guess).  The result, equivalent
to (\ref{vierie}) is:
\begin{equation}\label{eq:asymp1}
\log {v_n \over (2n)!} = 2n -4 \log n + C + {11\over 6 n} + {141\over
  160 n^2} + \dots,  
\end{equation}
where $C:= -3 -\log \pi -{3 \over 2} \log {8\over 3}$. In the
appendix we present a proof by Don Zagier of this asymptotic formula.

\section{Comparison with two other sequences}

As a matter of curiosity, we now compare our results so far with
similar results from two other sequences of enumerative geometry.
We shall see that the first case has quite similar features, while the
second case is more intricate.

\subsection{Numbers of plane rational curves}

One sequence of integers from enumerative geometry is $n_d$, the
number of plane rational curves of degree $d$ through $3d-1$ points in
$\P^2$.  Kontsevich's recursion formula \cite{KM-94} reads
$$
n_d =\sum_{k=1}^{d-1} n_k ~n_{d-k} \Big[
k^2 (d-k)^2 \left({3d-4\atop 3k-2}\right) -k^3 (d-k) \left({3d-4\atop 3k-1}\right)
\Big], \qquad\qquad n_1=1.
$$
The result is $n_1=1$, $n_2=1$, $n_3=12$, etc, ie. there is 1 line
through 2 points of the plane, 1 conic through 5 points of the plane,
12 cubics through 8 points of the plane, etc. 

We can similarly draw tables of $n_d$ mod $k$ for any integer $k$.
The results (in the same convention as before) are:
\begin{itemize} \setlength{\itemsep}{-.1cm} \setlength{\itemindent}{0cm}
\item[*] $k=2$: both rows vanish (except first two values), ie. all
  $n_d$ are even.
\item[*] $k=2^l$: all rows are 0, ie $n_d\equiv 0$ mod $2^l$ for $n>l+1$.
\item[*] $k=3$:  $n_{3d}\equiv 0$ mod 3, $n_{3d+2}\equiv 1$ mod 3,
  $n_{3d+1}\equiv $ alternating 1 or 2 mod 3 because $n_{6d+2}\equiv 4$ mod~6.
\item[*] $k=5$:  $n_{d}\equiv 0$ mod 5, for $d>8$.  Idem for $k=25$ ($d>23$)
\end{itemize}
Because of this big symmetry for low primes, most non-primes will
yield constant or regular rows (ie rows repeating when shifting
horizontally).  The only non-obvious case is $k=26$, where there is a
shift by 8 (because $k=13$ shifts by 16) and rows 4,6 alternate with 0.

Further, we only found three primes with regular features:
\begin{itemize} \setlength{\itemsep}{-.1cm} \setlength{\itemindent}{0cm}
\item[*] $k=7$: all rows are regular (repeat when shifted horizontally
  by 4), rows 5 and 7 are 0. 
\item[*] $k=13$: idem, shift by 16, no 0 row. 
\item[*] $k=19$: idem, shift by 12, no 0 row. 
\item[*] $k=5,11,17,23,29$: these primes give almost-0 rows
  (ie. $n_d\equiv 0$ except for a finite number of $d$).
\end{itemize}
We have not attempted to prove these observations.

\subsubsection{Asymptotics}
We now turn to the asymptotics of the sequence $n_d$ for $d\to\infty$. 
Di Francesco and Itzykson proved  \cite[Proposition 3]{FI} that
$$ {n_d\over (3d-1)!}\;=\;{A^d\over d^{7/2}}\,\biggl(B+{\rm O}\bigl({1\over d}\bigr)\biggr)\,, $$ 
as $d$ tends to infinity, and found the approximate values $A\approx 0.138$ and $B\approx 6.1$ for the
constants $A$~and~$B$.  Assuming a full asymptotic expansion 
$$ \frac{n_d}{(3d-1)!}\;\sim\;{A^d\over d^{7/2}}\biggl(B_0+\frac{B_1}d+\frac{B_2}{d^2}+\cdots\biggr)\,,$$
and applying variant~(II) of the asymp$_k$ trick to $\log (n_d/(3d-1)!)$, we obtain the much more accurate approximations
\begin{eqnarray*}  A\;&\approx&\; 0.138009346634518656829562628891755541716014121072\,, \\
                 B_0\;&\approx&\; 6.0358078488159024106383768720948935\,,\end{eqnarray*}
as well as the further values $B_1\approx-2.2352424409362074$,  $B_2\approx 0.054313787925$.
Unfortunately, we are not able to recognize any of these apparently irrational numbers, e.g., PARI does not see
in $\log A$ and $\log B_0$ a linear combination of simple numbers like 1, $\log 2$, $\log 3$, $\log\pi$, $\pi$ and $\pi^2$.

\subsection{Numbers of rational curves on the quintic threefold}

The other sequence we now introduce for the purpose of comparison is
$q_d$, that of holomorphic rational curves of degree $d$ embedded in
the quintic Calabi-Yau threefold.  These are the `instanton numbers'
of Candelas et al \cite{CdGP-91}. They are defined by the following line:
\begin{equation} \label{eq:def-instantons}
5+ \sum_{n\geq 1} q_n ~n^3 {q^n \over 1-q^n} = \left({q\over
    x}{dx\over dq}  \right)^3 ~{5\over (1-5^5 x) ~y_0(x)^2 } = 5 +2875
~q + 4876875 ~q^2 +\dots,  
\end{equation}
where $ q(x)= x~e^{\tilde{y}_1 /y_0} =x+770~x^2 +\dots $ is the
``mirror map'' and its inverse is $ x(q)= q -770~q^2 +\dots $. The
functions $y_0$ and  $\tilde y_1$ are solutions of a Picard-Fuchs
differential equation and are given by 
$$ 
y_0(x) := \sum_{n\geq 0} {(5n)! \over n!^5} ~x^n \qquad \textrm{and}
\qquad \tilde{y}_1(x)= \sum_{n\geq 0} \left( {(5n)! \over n!^5} ~5
\sum_{j=n+1}^{5n} {1\over j}  \right) ~x^n .
$$
When computing the numbers $q_d$, the longest step is without doubt
the inversion of the series $q(x)$.  The first few values are
$q_1=2875$, $q_2=609250$, etc.

We can again draw tables of $q_d$ mod $k$ for any integer $k$.
The main results (in the same convention as before) are:
\begin{itemize} \setlength{\itemsep}{-.1cm} \setlength{\itemindent}{0cm}
\item[*] $k=2$: second row is 0, ie. $q_{2d}$ are even.
\item[*] $k=4,8,16$: last row is 0, ie. $q_{2^l d} \equiv$ 0 mod $2^l$
  ($l\leq 4$) 
\item[*] $k=8,16$: rows 4,8,12,... are also 0, ie. $q_{2^l d+4m}
  \equiv$ 0 mod $2^l$ ($l\leq 4$)
\item[*] $k=32$: no 0 rows anymore!
\item[*] $k=5$: all rows are 0, idem at $k=25$.
\item[*] $k=20$: row 4,8,12,16,20 are 0, because both are 0 mod 4 and
  mod 5.
\end{itemize}
These congruences are much less impressive than in the previous cases,
due to the more complicated origin of the instantons. It is not even
mathematically understood why these are integers and what exactly they
count. Again we have not attempted to prove the congruences.

\subsubsection{Asymptotics}
In this case the asymptotic behaviour is much more tricky than in the
previous two examples. The growth is indeed exponential, but $\log
q_d$ has more than a simple logarithmic term and monomial terms.
Indeed, subtracting the logarithmic term gives us a sequence on which
the asymp$_k$ trick works badly -- as if other logarithmic terms were
hiding.  In fact, finding out the coefficient of the first log term is
already quite tough, and differentiating (to get rid of log-terms)
does not yield anything with only monomials.  
The second author's attempts to
deal with the asymptotics of $\log q_d$ through the asymp$_k$ trick can
be found in \cite{G-05}.

\section{Conclusion}

Though the congruences satisfied by the sequence $v_n$ of lines in
$\P^n$ are numerous and very interesting, the exponential asymptotic
growth $v_n \sim n^{2n}$ rules out that the $v_n$ are Fourier
coefficients of any modular form on a subgroup of $SL(2,\Z)$ (whose
coefficients typically grow like $n^{2k-1}$ for weight $2k$).  It is
quite gratifying to see that the asymptotic expansion in
(\ref{eq:asymp1}) can be written out with as many exact terms as one
wishes, since the coefficients are rational.

Another sequence -- that of numbers $n_d$ of degree $d$ curves through
$3d-1$ points of the plane -- has very similar behaviour, in terms of
congruences as well as asymptotics (though the latter's coefficients
will be irrational and not recognized).  A last sequence -- that of
numbers $q_d$ of degree $d$ rational curves on the quintic Calabi-Yau
threefold -- is much less enticing; its congruences are rather limited
and its asymptotics are awkward: not simply one or two log-terms
followed by mere monomial terms, but certainly log(log) terms or
infinitely many log terms.  

It would be interesting to study more of typical sequences from
enumerative geometry and see if there is an underlying pattern.  Also,
the question of how many of those sequences satisfy a recurrence
relation is still open.  Among the three sequences that we discussed,
only that of $n_d$ (plane rational curves) obeys a known
recurrence.

\subsection*{Acknowledgments}

The first author, who initiated this paper, 
is deeply grateful for Don Zagier's help throughout all the stages
of this project. Not only did Don simplify matters considerably  with his
formula (\ref{eq:defn}), but he also equipped him with the asymp$_k$
trick.  He is thankful
for fruitful discussions with Robert Osburn. The second author likes
to thank Daniel Berend, Paolo Dominici and Yossi Moshe for interesting e-mail correspondence 
and Carl Pomerance for providing us with Lemma \ref{carl}. He especially likes
to thank Don Zagier for his help in improving the exposition of the paper.
Dmitry Kerner and Masha Vlasenko we thank for several helpful discussions.
Several
inaccuracies in an earlier version were detected by Alexander Blessing during
a two week practicum he did with the second author.\\
\indent This project was made possible thanks to the support of the MPIM in
Bonn.\\

\vskip1cm
\centerline{\bf APPENDIX}
\centerline{\it by Don Zagier}
\centerline{\bf Exact and asymptotic formulas for $v_n$}
\label{sec:appendix}
\vskip.5cm
\noindent In this appendix we prove the alternative definition (\ref{eq:defn}) and the
asymptotic formulae (\ref{vierie}) and (\ref{eq:asymp1}) for the numbers $v_n$ defined in~(\ref{sommie}).
\vskip.5cm

\centerline{\bf Exact formulas}\medskip

\begin{prop}
  Let $G(x,y)$ be a homogeneous polynomial of degree $2n$ in two variables and $P(x)$ a monic
  polynomial of degree $n+1$ with distinct roots.  Then the expression
\begin{equation}
\label{A1}
\sum_{\alpha,\,\beta\in\C\atop P(\alpha) =P(\beta)=0}
  \frac{G(\alpha,\beta)}{P'(\alpha)P'(\beta)} 
\end{equation}
is independent of $P$ and equals the coefficient of $x^ny^n$ in $G(x,y)$.  
\end{prop}
{\it Proof}. By linearity it is enough to consider monomials $G(x,y)=x^r y^s$, $r+s=2n$.
  Then the expression~(\ref{A1}) factors as
  $\bigl(\sum_{P(\alpha)=0}\frac{\alpha^r}{P'(\alpha)}\bigr)
  \bigl(\sum_{P(\beta)=0}\frac{\beta^s}{P'(\beta)}\bigr)$.  But by the
  residue theorem we have
$$
\sum_{P(\alpha)=0}\frac{\alpha^r}{P'(\alpha)} =\sum_{\alpha\in\C}
\textrm{Res}_{x=\alpha}\biggl(\frac{x^r\,dx}   
{P(x)}\biggr)=-\textrm{Res}_{x=\infty}\biggl(\frac{x^r\,dx}{P(x)}\biggr)\,,
$$
and this equals $0$ if $0\le r<n$ and $1$ if $r=n$ since $P$ is monic of
degree $n+1$. The proposition follows. 
\qed\\

\noindent {\tt Remark.}  The same proof shows that if $G$ is homogeneous of
degree $m+n$ and $P$ and $Q$ are two monic polynomials of degrees $m+1$ and $n+1$
with distinct roots, then $$\sum_{P(\alpha) =Q(\beta) =0}
\frac{G(\alpha,\beta)} {P'(\alpha) Q'(\beta)}$$ is independent of $P$
and $Q$ and is equal to the coefficient of $x^my^n$ in $G(x,y)$.  Yet
more generally, and still with the same proof, if $G$ is a homogeneous polynomial
of degree $n_1+\cdots+n_k$ 
in $k$ variables and $P_1,\dots,P_k$ monic polynomials of degree
$n_1+1, \dots, n_k+1$ with no multiple roots, then
$$\sum_{P_1(\alpha_1) =\dots= P_k(\alpha_k)=0} \frac{G(\alpha_1,\dots,
  \alpha_k)} {P_1'(\alpha_1) \cdots P_k'(\alpha_k) }$$ is independent
of all the $P_i$ and is equal to the coefficient of $x_1^{n_1 }\cdots
x_k^{n_k}$ in $G(x_1,\dots, x_k)$. In fact $G$ need not even be
homogeneous, but can be any polynomial in $k$ variables of degree
$\le n_1+\cdots+n_k$.

\begin{cor}
  Let $F(x,y)$ be a symmetric homogeneous polynomial of degree $2n-2$
  in two variables and $w_0,\ldots,w_n$ distinct complex numbers.
  Then the expression
  $$
  \sum_{0\le i<j\le n}{F(w_i,w_j)\over {\displaystyle{\prod_{0\le k\le
    n\atop k\ne i,\,j}(w_i-w_k)(w_j-w_k)}}}
  $$
  is independent of $w_0,\ldots,w_n$ and equals the coefficient of
  $x^{n-1}$ in $(1-x)F(x,1)$.
\end{cor}
{\it Proof}. This follows after a short
calculation if we apply the proposition to $G(x,y)=(x-y)^2F(x,y)$, $P(x)=\prod_{i=0}^n(x-w_i)$. \qed\\

Corollary 2 immediately implies that the right hand side of equation (\ref{sommie}) is independent
of the (distinct) complex variables $w_0,\ldots,w_n$ and that (\ref{sommie}) is equivalent to (\ref{eq:defn}). The
computational advantage is huge: formula (\ref{sommie}) is very slow to compute, even for
moderately large $n$, whereas (\ref{eq:defn}) can be implemented in
PARI in one line as
$$
\textrm{\tt v(n) = coeff(prod(j=0,2*n-3,2*n-3-j+j*x,1-x),n-1) }
$$
and takes $<2$ seconds to compute $v_n$ up to $n=100$ and 46
seconds up to $n=224$.

We can rewrite (\ref{eq:defn}) in several other forms by using residue calculus. Setting $D=2n-3$
and making the substitution $x=1-D/z$, we find
\begin{equation}
\label{A3}
~~~~~~~~~v_n=\textrm{Res}_{x=0}\biggl((1-x)\, \prod_{j=0}^{2n-3}(2n-3-j+jx)\;
 \frac{dx}{x^n}\biggr)
\end{equation}
\begin{equation}
\label{A4}
=D^{2n}\;\textrm{Res}_{z=D}\biggl(\frac{\prod_{j=0}^D(z-j)}{z^{n+1}\,(z-D)^n}\,dz\biggr).
\end{equation}
Since the residue of the integrand at infinity is zero, we can also
write this as  
\begin{equation}
\label{A5}
v_n =-\,D^{2n}\; \textrm{Res}_{z=0}\biggl( \frac{\prod_{j=0}^D(z-j)}
 {z^{n+1}\, (z-D)^n}\,dz  \biggr),
\end{equation}
while simply making the substitution $z\mapsto D-z$ in (\ref{A4})
gives the similar expression
\begin{equation}
\label{A6}
v_n = D^{2n} \textrm{Res}_{z=0} \biggl( \frac{\prod_{j=0}^D(z-j)}
 {z^{n} ~(z-D)^{n+1}}~ dz  \biggr),
\end{equation}
and adding these two last expressions gives yet a third form:
\begin{equation} 
\label{A7}
v_n = \half D^{2n+1} \textrm{Res}_{z=0} \biggl( \frac{\prod_{j=0}^D(z-j)}
 {z^{n+1} ~(z-D)^{n+1}}~ dz  \biggr).
\end{equation}
Each of the formulae (\ref{A5})--(\ref{A7}) expresses $v_n$ as
the constant term at $z=0$ of the Laurent expansion of a rational
function, e.g.~(\ref{A5}) says
\begin{equation}
\label{A8}
v_n =(-1)^n D^{2n}\cdot\textrm{coefficient of $z^{n-1}$ in
 $\frac{(1-z)(2-z) \cdots(D-1-z)}{(D-z)^{n-1}}$ as $z\to 0$.} 
\end{equation}
Substituting $z=Du$, we can write this as
\begin{equation}
\label{A9}
v_n=(-1)^nD^2\cdot\textrm{coefficient of $u^{n-1}$ in
 $\frac{(1-Du)(2-Du) \cdots(D-1-Du)}{(1-u)^{n-1}}$ as $u\to 0$,}
\end{equation} 
 from which we see again that $D^2|v_n$ (Lemma \ref{lemma0}).
By expanding $(D-z)^{1-n}$ by the binomial theorem, we also obtain closed formulae
for $v_n$; for instance (\ref{A8}) gives:
\begin{equation}
\label{A10}
  v_n =\sum_{m=0}^{n-1} (-1)^{n-1-m} {2n-2-m \choose n-1} D^{m+1}
  \Bigg[{D\atop m}\Bigg],
\end{equation}
where $\Big[{D\atop m}\Big]$, the coefficient of $z^m$ in
$z(z+1)\cdots (z+D-1)$, is a Stirling number of the first kind.\\

\centerline{\bf Asymptotics}
\noindent To obtain the asymptotic expansion of $v_n$, we write the residue in
(\ref{A3}) as $\frac{1}{2\pi i} \int_{|x|=1}$ and we make the
substitution $x=(1+it)/(1-it)$ to obtain, after a short calculation,
\begin{equation}
\label{A11}
v_n =\frac{2}{\pi} \int_{-\infty}^{\infty} \prod_{r=1,3,\dots,D} \! \Big(
\frac{D^2 +r^2 t^2}{1+t^2} \Big) \, \frac{t^2 ~dt}{(1+t^2)^2}
=\frac{2}{\pi} D^{D+1}  \int_{-\infty}^{\infty} \phi_D(t) \frac{t^2
  ~dt}{(1+t^2)^2},
\end{equation} 
where $\phi_D(t)$ denotes the rational function
$$
\phi_D(t) = \prod_{r=1,3,\dots,D} \frac{1+r^2 D^{-2} t^2}{1+t^2}.
$$
It is easy to see that $\phi_D(0)=1$ and $\phi_D(t)\le e^{-cDt^2}$ for some absolute constant $c>0$ (a much
more precise formula will be given in a moment), so the main contribution to the integral comes
from small $t$.  For $t$ small and $D$ large we have (uniformly in both variables)
\begin{eqnarray} 
&&\log \phi_D(t) = \sum_{j=1}^\infty\frac{(-1)^{j-1}}j\,\biggl[\,\sum_{r=1,3,\ldots,D}
  \biggl(\frac{r^{2j}}{D^{2j}}-1\biggr)\biggl]\,t^{2j} \nonumber\\
&&= \biggl(-\frac D3+\frac1{3D}\biggr)\,t^2+\biggl(\frac D5-\frac1{3D}+\frac2{15D^3}\biggr)\,t^4
+\, \biggl(-\frac D7+\frac1{3D}-\frac4{9D^3}+O\biggl(\frac1{D^5}\biggr)\biggr)\,t^6\nonumber\\
&&+\, \biggl(\frac D9-\frac1{3D}+\frac{14}{15D^3}+O\biggl(\frac1{D^5}\biggr)\biggr)\,t^8
+\,\biggl(-\frac D{11}+\frac1{3D}+O\biggl(\frac1{D^3}\biggr)\biggr)\,t^{10}\nonumber\\
&&+\, \biggl(\frac D{13}-\frac1{3D}+O\biggl(\frac1{D^3}\biggr)\biggr)\,t^{12}
+\, \biggl(-\frac D{15}+O\biggl(\frac1D\biggr)\biggr)\,t^{14}
+\, \biggl(\frac D{17}+O\biggl(\frac1D\biggr)\biggr)\,t^{16}+O\bigl(Dt^{18}\bigr),\nonumber
\end{eqnarray}
and hence
\begin{eqnarray} 
&\,&\frac{x^2}{(1+x^2/D)^2}\, \phi _D\bigl(\frac x{\sqrt D}\bigr)\nonumber\\
&\,&= e^{-x^2/3}\biggl[x^2+\biggl(\frac{x^6}5-2x^4\biggr)D^{-1}
      +\biggl(\frac{x^{10}}{50}-\frac{19x^8}{35}+3x^6+\frac{x^4}3\biggr)D^{-2}  \nonumber\\
&\,&+\, \biggl(\frac{x^{14}}{750}-\frac{12x^{12}}{175}
      +\frac{314x^{10}}{315}-\frac{59x^8}{15}-x^6\biggr)D^{-3}+\cdots \nonumber\\
&\,&+\, \biggl(\frac{x^{30}}{393750000}-\frac{11x^{28}}{19687500}+\cdots+\frac{355x^{10}}{162}
      +\frac{2x^8}{45}\biggr)D^{-7}+O\bigl(D^{-8}\bigr)\biggr].\nonumber
\end{eqnarray}
Substituting this expansion (with the 34 omitted terms included) into equation (\ref{A11}) with $t$
replaced by $x/\sqrt{D}$ and using the standard evaluation
$$\int_{-\infty}^\infty e^{-x^2/3} x^{2n} dx =\frac{(2n)!}{n!} \Big(\frac{3}{4}\Big)^n \sqrt{3\pi},$$ we obtain
\begin{eqnarray}
v_n&=&\sqrt{\frac{27}\pi} D^{D-1/2} \biggl(1-\frac94 D^{-1}+\frac{969}{160}D^{-2}
  -\frac{61479}{3200}D^{-3} + \frac{25225773}{358400}D^{-4} \nonumber\\
  & -& \frac{10092025737}{35840000}D^{-5}+\, \frac{2271842858513}{2007040000}D^{-6}
  -\, \frac{4442983688169}{1146880000}D^{-7}+O\bigl(D^{-8}\bigr)\biggr).\nonumber
\end{eqnarray}
This asymptotic formula can of course be written in many other ways,
e.g.: 
$$
v_n=  \sqrt{\frac{27}{\pi}} (2n-3)^{2n-7/2} \Big( 1-\frac{9}{8n}
-\frac{111}{640 n^2} -\frac{9999}{25600 n^3}+\frac{87261}{5734400n^4} -\cdots\Big)
$$
or
$$
v_n =e^{-3}  \sqrt{\frac{27}{\pi}} (2n)^{2n-7/2} \Big( 1+\frac{15}{8n}
+\frac{1689}{640 n^2} +\frac{79281}{25600 n^3} +\frac{19691853}{5734400n^4}+\cdots\Big)
$$
or
$$
\log \frac{v_n}{(2n)!} =2n-4\log n +C +\frac{11}{6n} +\frac{141}{160
  n^2} +\frac{9973}{28800 n^3} +\frac{59673}{179200n^4}+\cdots 
$$
with $C=-3-\log \pi -\frac{3}{2}\log\frac{8}{3}$. Of course, more
terms could be obtained in any of these expansions if desired.

\medskip\noindent {\footnotesize 
{\tt D.~Gr\"unberg}, 49 rue Fondary,
75015 Paris,
e-mail: {\tt grunberg@mccme.ru}\\
{\tt P.~Moree}, Max-Planck-Institut f\"ur Mathematik,
Vivatsgasse 7, D-53111 Bonn, Germany.\\ e-mail: {\tt moree@mpim-bonn.mpg.de}\\
{\tt D.~Zagier}, Max-Planck-Institut f\"ur Mathematik,
Vivatsgasse 7, D-53111 Bonn, Germany.\\ e-mail: {\tt don@mpim-bonn.mpg.de}}


{\small
\begin{thebibliography}{99}

\bibitem{Br} F.L.~Bauer, For all primes greater than $3,\;({2p-1\atop p-1})\equiv 1 \;({\rm mod~}p^3)$ holds,
{\it Math. Intelligencer} {\bf 10} (1988), no. 3, 42.

\bibitem{BH} D.~Berend and J.E.~Harmse, On some arithmetical properties of middle binomial coefficients, 
{\it Acta Arith.} {\bf 84} (1998), 31--41.

\bibitem{Beu} F.~Beukers, Irrationality proofs using modular forms, Ast\'erisque No. {\bf 147}-{\bf 148} (1987), 271--283.


\bibitem{Cai}  T.~Cai, A congruence involving the quotients of Euler and its applications. I, 
{\it Acta Arith.}  {\bf 103}  (2002),  313--320.   

\bibitem{CdGP-91} P.~Candelas, J.X.~de la Ossa, P.~Green,
  and L.~Parkes, {\it A Pair of Calabi-Yau Manifolds as an Exactly Soluble
  Superconformal Field Theory},  Nucl. Phys. {\bf B359} (1991), 21--74,
  reprinted in {\it Essays on Mirror Manifolds} (S.T.~Yau, ed.)
  Hong Kong, 1992.
  
\bibitem{FI} P.~Di Francesco and C.~Itzykson, Quantum intersection rings.  {\it The moduli space of curves} 
(Texel Island, 1994),  81--148, Progr. Math., 129, Birkh\"auser Boston, Boston, MA, 1995.   
  

\bibitem{Dom} P.~Dominici, Sequence A027363, {\it On-Line Encyclopedia of Integer Sequences},
{\tt http://www.research.att.com/\~{}njas/sequences/}.  


\bibitem{F-84} W.~Fulton,  {\it Intersection Theory}, Springer, New York, 1984.
 
\bibitem{AG} A.~Granville, Arithmetic properties of binomial coefficients. I. Binomial coefficients modulo prime powers, 
{\it Organic mathematics} (Burnaby, BC, 1995),  253--276, CMS Conf. Proc., 20, Amer. Math. Soc., Providence, RI, 1997. 
 
 \bibitem{G-04}  D.~Gr\"unberg, Integrality of open instanton
  numbers, {\it J. Geom. Phys.} {\bf 52} (2004), 284--297, hep-th/0305057.
  
\bibitem{G-05}  D.~Gr\"unberg, {\it Asymptotic growth of the sequence of instantons on
the quintic threefold}, unpublished manuscript, 2005.

\bibitem{HW} G.H.~Hardy and E.M.~Wright, {\it An introduction to the theory of numbers}. Fifth edition. The Clarendon 
Press, Oxford University Press, New York, 1979.
  
\bibitem{KS-92} S.~Katz and S.A.~Stromme,  {\it Schubert -a Maple
 package for intersection theory and enumerative geometry}, (1992)
 {\tt http://www.mi.uib.no/\~{}stromme/schubert/}

\bibitem{KM-94} M.~Kontsevich and Y.~Manin, Gromov-Witten
    classes, quantum cohomology, and enumerative geometry,
    {\it Comm. Math. Phys.} {\bf 164} (1994), 525--562, hep-th/9402147.
    
\bibitem{Mani} L.~Manivel, 
{\it Symmetric functions, Schubert polynomials and degeneracy loci}, 
SMF/AMS Texts and Monographs {\bf 6}, American Mathematical Society, Providence, RI, 2001.

\bibitem{M1} Y.~Moshe, The density of 0's in recurrence double sequences, 
{\it J. Number Theory}  {\bf 103}  (2003),  109--121.

\bibitem{M2} Y.~Moshe, The distribution of elements in automatic double sequences, 
{\it Discrete Math.} {\bf 297}  (2005),  91--103.



\bibitem{SB} J.~Stienstra and F.~Beukers, On the Picard-Fuchs equation and the formal Brauer group of certain 
elliptic $K3$-surfaces, {\it Math. Ann.} {\bf 271} (1985), 269--304.

    
\bibitem{vdW0} B.L.~van der Waerden, Zur algebraischen Geometrie. II. 
Die geraden Linien auf den Hyperfl\"achen des ${\mathbb P}_n$.
{\it Math. Ann.} {\bf 108}  (1933),  253--259.
    
\bibitem{vdWS} B.L.~van der Waerden, {\it Zur algebraischen Geometrie. 
Selected papers}, Springer-Verlag, Berlin, 1983.

\bibitem{Z} D.~Zagier, Vassiliev invariants and a strange identity related to
the Dedekind eta-function, {\it Topology} {\bf 40} (2001), 945--960.





\end{thebibliography}}
\end{document}